\documentclass[a4paper,12pt]{article}

\usepackage[T1]{fontenc}
\usepackage[utf8]{inputenc}
\usepackage[english]{babel}
\usepackage[inline,shortlabels]{enumitem}
\usepackage{amsmath}
\usepackage{amssymb}
\usepackage{amsthm}
\usepackage[noadjust]{cite}
\usepackage[a4paper,total={6in,9in}]{geometry}
\usepackage[colorlinks,citecolor=blue]{hyperref}
\usepackage[capitalize,noabbrev,nameinlink]{cleveref}
\usepackage{xcolor}
\usepackage{longtable}
\usepackage{tikz,pgfplots,pgfplotstable}
\usetikzlibrary{patterns}
\usepackage{algorithmic}
\usepackage{marginnote}
\usepackage[labelformat=simple]{subcaption}
\usepackage{multirow}
\usepackage[tableposition=above]{caption}
\setlist{font=\normalfont,topsep=1ex,parsep=0ex}

\setlist[enumerate]{label=(\alph*)}
\newcommand{\smallsection}[1]{%
  \par\medskip
  {\normalsize\bfseries #1\par\nobreak\smallskip}
  \normalsize
}

\numberwithin{equation}{section}
\numberwithin{table}{section}    
\numberwithin{figure}{section}

\crefname{figure}{Figure}{Figures}
\crefname{table}{Table}{Tables}
\crefname{assumption}{Assumption}{Assumptions}
\Crefname{ALC@unique}{Step}{Steps}

\newlist{alglist}{enumerate}{1}
\setlist[alglist]{topsep=1ex,parsep=0ex,leftmargin=*,label=\textbf{Step~\arabic*.}}

\newcommand{\R}{\mathbb{R}}

\newcommand\norm[1]{\left\Vert#1\right\Vert}

\newcommand{\toattentive}[1]{\overset{#1}{\to}}

\DeclareMathOperator{\dist}{dist}

\newcommand{\dom}{\operatorname{dom}}
\newcommand{\identity}{\operatorname{Id}}

\newtheoremstyle{bolddef}{}{}{\normalfont}{}{\bfseries}{.}{ }{\thmname{#1}\thmnumber{ #2}\thmnote{ (#3)}}

\newtheoremstyle{boldplain}{}{}{\itshape}{}{\bfseries}{.}{ }{\thmname{#1}\thmnumber{ #2}\thmnote{ (#3)}}

\theoremstyle{bolddef}
\newtheorem{definition}{Definition}[section]
\newtheorem{algorithm}[definition]{Algorithm}
\newtheorem{assumption}[definition]{Assumption}

\theoremstyle{boldplain}
\newtheorem{lemma}[definition]{Lemma}
\newtheorem{theorem}[definition]{Theorem}
\newtheorem{proposition}[definition]{Proposition}

\newlength\figureheight
\newlength\figurewidth

\pgfplotsset{width=7cm,compat=1.3}

\usepackage{algorithm}

\definecolor{todocolor}{rgb}{1.0,0.0,0.0}

\newcommand\email[1]{\href{mailto:#1}{\texttt{#1}}}

\newcommand{\orcid}[1]{ORCID: \href{https://orcid.org/#1}{#1}}

\newcommand{\mscLink}[1]{\href{http://www.ams.org/mathscinet/msc/msc2020.html?t=#1}{#1}}

\begin{document}

\title{
	\bfseries\scshape 
	General Proximal Quasi-Newton Methods based on model functions for nonsmooth nonconvex problems
	}

\author{Xiaoxi Jia%
	\thanks{%
	Independent researcher,
		\email{xiaoxijia26@163.com},
		\orcid{0000-0002-7134-2169}
	}
}

\maketitle
{
\small\textbf{\abstractname.}
In this manuscript, we propose a general proximal quasi-Newton method tailored for nonconvex and nonsmooth optimization problems, where we do not require the sequence of the variable metric (or Hessian approximation) to be uniformly bounded as a prerequisite, instead,  the variable metric is updated by a continuous matrix generator.  From the respective of the algorithm, the objective function is approximated by the so-called local model function and subproblems aim to exploit the proximal point(s) of such model function, which help to achieve the sufficiently decreasing functional sequence along with the backtracking line search principle.  Under mild assumptions in terms of the first-order information of the model function,  every accumulation point of the generated sequence is stationary and the sequence of the variable metric is proved not to be bounded.  Additionally, if the function has the Kurdyka-{\L}ojasiewicz property at the corresponding accumulation point,  we find that the whole sequence is convergent to the stationary point, and the sequence of the variable metric is proved to be uniformly bounded.  Through the above results, we think that the boundedness of the sequence of the variable metric should depend on the regularity of objectives,  rather than being assumed as a prior for nonsmooth optimization problems.  Numerical experiments on polytope feasibility problems and (sparse) quadratic inverse problems demonstrate the effectiveness of our proposed model-based proximal quasi-Newton method,  in comparison with the associated model-based proximal gradient method.
\par\addvspace{\baselineskip}
}

{
\par\addvspace{\baselineskip}
}

{
\par\addvspace{\baselineskip}
}

\section{Introduction}
\indent Let us consider
\begin{equation}\label{Eq:P}\tag{P}
	\min_{x\in \mathbb R^n} \ f(x),
\end{equation}
where $f: \mathbb R^n \to \overline{\mathbb R}$ is assumed to be proper and lower semicontinuous.  A natural approach is to minimize the approximation of $f(x)$,  which has better structures than $f$, such approximation is commonly referred to as the so-called model function.  The essence of utilizing model functions lies in controlling their distance from the actual objective function,  striving for it to be sufficiently small.  In fact,  we just know the information about the function value $f$, not about $\nabla f$ at al.  Then we settle for less and use the subdifferential of $f$.  In particular, when $f$ is smooth, the most prevalent model example is the first-order Taylor approximation. In nonsmooth optimization, the key consideration is the approximation quality of the model function (or approximation error), which is usually controlled by the so-called (real-valued) growth function to quantify the approximation error at the current iterate.  Centered at some $\bar x \in \dom f$,  it can be formulated mathematically as
\begin{equation}\label{Eq:M}
|f_{\bar x}(x)-f(x)| \leq \omega(\|x-\bar x\|) \quad \forall x\in \dom f,
\end{equation} 
where $f_{\bar x}: \mathbb R^n \to \overline{ \mathbb R}$ is the model function centered at $\bar x$ and $\omega: \mathbb R_+ \to \mathbb R_+$ is the growth function. Drusvyatskiy et al. \cite{drusvyatskiy2021nonsmooth} characterized such model functions $f_{\bar x}: \mathbb R^n \to \overline{\mathbb R}$ as Taylor-like models.  In this manuscript,  we introduce a revised definition of the model function (see \Cref{Def:GrowthFun}), where we relax the model approximation principle \eqref{Eq:M} merely for those $x$ in some neighborhood of $\bar x$, not for all $x \in \dom f$.  This offers a more flexible framework to better capture characteristics of the original functions, in particular, which exhibit the nonsmooth behaviour or whose  gradients are just locally Lipschitz continuous.  

Note that the exact minimization of the model function is  possibly ineffective.  Instead, the model function is typically complemented by a proximity measure, which encourages solutions closed to the current iterate.  Consequently,  subproblems arise, wherein the objective function  becomes the sum of the model function and the proximity measure. When $f$ is smooth,  the model function can be read using its gradient information,  we then employ the Euclidean norm as a proximity measure and transform computing the next iterate into a gradient decent step,  and for the norm deduced by some variable metric as a proximity measure,  we normally employ the proximal (quasi-)Newton methods to solve the actual problem.  Bregman in \cite{bregman1967relaxation} proposed to invoke a more general proximity measure as afforded by the so-called Bregman distances.  Based on this idea,  minimization of subproblems results in Bregman proximal algorithm's update step \cite{BolteSabachTeboulleVaisbourd2018,ochs2019non,mukkamala2022global}.  However, the choice of the Bregman function is problem-dependent and non-trivial,  because this significantly impacts the efficiency of subproblems which sometimes require pretty complicated solvers.  In comparison, the proximity deduced by the variable metric are simultaneously powerful and simple.  


This manuscript focuses on proximal quasi-Newton methods for solving \eqref{Eq:P}, where, in each step, the subproblem
\begin{equation}\label{Eq:introsubp}
\min_{x} f_{x^k}(x)+\frac{\gamma_k}{2}(x-x^k)^TH_k(x-x^k),
\end{equation}
where $x^k$ denotes the current iterate and $1/\gamma_k$ is the stepsize,  needs to be solved. In order to ensure the global convergence, we integrate the solution(s) of the subproblem with a backtracking line search technique.  Note that the choice of the matrix $H_k$ is crucial for developing such algorithms. For example, first-order methods use $H_k$ as a positive multiple of the identity matrix,  while (quasi-)Newton methods denote $H_k$ as the (approximated) Hessian.   In \cite{mukkamala2022global}, Mukkamala et al.  proposed a classical proximal gradient method ($H_k:=\identity$), where the boundedness of the corresponding iterative sequence is needed for the convergence analysis,  which can be achieved by essentially requiring that $f$ is coercive. 
Some proximal Newton methods or proximal quasi-Newton methods also require the sequence $\big(x^k\big)_{k \in \mathbb N}$ to be bounded for the desired convergence results, cf. \cite{mordukhovich2023globally, kanzow2022efficient}.  In this manuscript,  let us now emphasize that we do not assume the boundedness of $\big(x^k\big)_{k \in \mathbb N}$ at all.  Additionally, many works involving the convergence and the rate-of-convergence of Newton-type methods normally rely on the regularity of $f$.  In particular,  requiring  that $f$ is (partially) smooth, and that the associated gradient to be Lipschitz continuous serves for the desired Q-linear convergence in 
 \cite{mordukhovich2023globally,lee2014proximal,chen1999proximal,scheinberg2016practical},  and sometimes the associated Hessian is Lipschitz continuous \cite{mordukhovich2023globally,lee2014proximal,nesterov2006cubic} for the desired Q-superlinear even Q-quadratic convergence rate.  However,  in our case, $f$ is merely assumed to be lower semicontinuous,  consequently we can not exploit the gradient and also the second-order information,  leading to more complicated (even potentially failed) convergence analysis for the proposed algorithm, particularly regarding the rate of convergence.  To overcome such challenges, we give a mild assumption about the first-order information,  then the subsequential convergence will be obtained, where the sequence of Hessian approximations is not bounded.  In order to derive the whole sequential convergence, we employ the Kurdyka–Łojasiewicz (KL) property \cite{Lojasiewicz1963, Lojasiewicz1965}, which naturally holds if the potential function is semialgebraic \cite{AttouchBolteRedontSoubeyran2010}.  After the KL property, the sequence of variable metrics is proved to be uniformly bounded.

Note that we initially proposed the simplest case of unconstrained minimization, which serves as the foundation for our broader programs encompassing various practical and interesting problems, such as the (addictive) composite problems \cite{lewis2016proximal,KanzowMehlitz2022,jia2023convergence}, difference of convexity \cite{Francisco,kanzow2023bundle},  fractional optimization problems \cite{crouzeix1985algorithm,dinkelbach1967nonlinear} and so on.  Depending on the choice of the approximate Hessian and the model function at iterations, our proposed algorithm (\Cref{Sec:Alg and Res}) covers many classical (sub)gradient methods \cite{Beck2017} and second-order methods \cite{nocedal1999numerical}. Importantly we do not impose any convex assumption on the objective function, making our work in this manuscript more general.

\section{Contributions}\label{Sec:Contr}
\textbf{\emph{Local} model function.}
As previously mentioned, for smooth functions, the model function is always chosen as the Taylor's approximation, which is unique.  For nonsmooth functions, there are only ``Taylor-like” model functions \cite{drusvyatskiy2021nonsmooth}.  Convex model functions are explored by Ochs et al.  in \cite{ochs2019non} and Ochs and Malitsky in \cite{ochs2019model}. Nonconvex model functions are discussed by Mukkamala in \cite{mukkamala2022global} and by Drusvyatskiy in \cite{drusvyatskiy2021nonsmooth}.  In the nonconvex case, previous work required a global control on the model approximation error by the so-called growth function.  
This concept is a generalization of the (global) Lipschitz or Hölder continuity \cite{ochs2019non}.  However, for functions without the global uniform continuity property,  one might fail to find the corresponding model function. In other words, the model approximation error, sometimes, cannot be captured by a globally uniform growth function. 

For example,  regarding a continuously differentiable function,  its Taylor-like model is always popular.  However,  even for such function with local Lipschitz continuous gradients,  
 the descent lemma yields
\begin{equation*}
\left|f(x^k)+\left<\nabla f(x^k), x-x^k \right>-f(x)\right| \leq \frac{L_{\bar x, x}}{2}\|x-\bar x\|^2,
\end{equation*}  
where $L_{\bar x, x}$ is the (local) Lipschitz parameter dependent on $\bar x$ and $x$. Since $\sup_{x\in \dom f}L_{\bar x, x}$ might be infinite,  we possibly fails to find a growth function such that the corresponding model approximation error can be bounded for all $x$ in the entire domain, although the used first-order model function is very classical.  If we now urge $x$ in some neighborhood, then the corresponding growth function dependent on the neighborhood center always exists.  This motivates us to relax the classical definition of the model function into its local version, please see \Cref{Def:ModelFunc}, which suits for much broader functions.  Precisely, let us take the function $x^4 \ (x\in \mathbb R)$  as an example,  whose gradient is Lipschitzly continuous.  Its first-order Taylor approximation centered at $\bar x \in \mathbb R$ involves a term $\bar x^3$,  yielding a local approximation that does not work globally.  Meanwhile, due to the local approximation principle, we can choose its model function as
\begin{equation*}
f_{\bar x}(x):=\max\{0, \bar x^4+4{\bar x^3}(x-\bar x)\},
\end{equation*}
which generates a better approximation than the classical first-order Taylor expansion.  Consider a composite problem $f(x):=|x^4-1|$, the Lipschitz continuity of the gradient is invalid. But, the local approximation principle holds if we choose the model function as
\begin{equation*}
f_{\bar x}(x):=|\bar x^4-1+4{\bar x^3}(x-\bar x)|.
\end{equation*}
Another model function with better structures can be given by
\begin{equation*}
f_{\bar x}(x):=\max\{0, \bar x^4-1+4{\bar x^3}(x-\bar x)\},
\end{equation*}
which also ensures the local approximation principle valid.

On the other hand, we are the first to generalize the subdifferential relationship between the model function and its original function (\cref{Prop: IRelation of subd}), which provides very vital informations for the convergence analysis in algorithms-driven situations.

\textbf{\emph{Unbounded} variable metric a priori.} Another significant contribution lies in the relaxation of the requirement for the variable matric to be upper bounded when employing quasi-Newton methods to solve \eqref{Eq:P}. Traditionally,  the works on quasi-Newton methods often assume the boundedness of the variable metric \cite{mordukhovich2023globally,scheinberg2016practical,lee2014proximal,kanzow2023regularization,kanzow2022efficient}.  But,  for nonsmooth problems with differentiabilities,  in particular, the objective function is locally Lipschitz continuous, the difference of the gradients might be enormous compared to the difference of the iterates, the inverse Hessian approximation typically becomes very ill-conditioned. Its eigenvectors corresponding to tiny eigenvalues are directions along which the function varies nonsmoothly \cite{lewis2013nonsmooth}.  Hence,  in nonsmooth optimization, the bounded assumption on the Hessian approximation results in the ineffectiveness even failure of quasi-Newton methods.  In some sense, the following references can confirm the hypothesis in nonsmooth optimization, they are,  Leconte and Orban in \cite{leconte2023complexity} noted that ``In the present paper,  we examine the situation where the sequence of Hessian approximations is allowed to grow unbounded”.  \cite[Section~8.4]{conn2000trust} proposed that for BFGS and SR1 approximations,  the Hessian approximation could potentially grow at a constant at each update, though it remains not clear whether that bound is achieved.  In practice, this assumption is pretty restrictive, consider the simple example of $x^a$, with $0<a<2$, $x\neq 0$.  Clearly, whenever $x \to 0$, we have that $| f^{''}(x)|=|a(a-1)x^{a-2}| \to \infty$.  

\textbf{Contributions when using the KL property.} In this manuscript, we do not assume the boundedness of the iterates, which is vital for the technical proof when using the KL property and has been required by many relevant publications \cite{BolteSabachTeboulleVaisbourd2018,ochs2019unifying,CohenHallakTeboulle2022,wu2021inertial}.  
In addition, 
we just require the sufficiently decreasing functional sequence and the local relative error condition, do not require the continuity condition like \cite{CohenHallakTeboulle2022,ochs2019unifying}.  
In the situation where the sequence of variable metric is unbounded, we demonstrate that the (whole) sequence generated by the Newton-like methods is convergent to a stationary point when employing the KL property.  Subsequently, our findings challenge the assertion made in \cite{stella2017forward} that the boundedness of the Hessian approximation is a prerequisite for the sequential convergence when the KL property is used.

In this manuscript,  we require the variable metric to be generated by a  continuous generator,  do not assume the sequence of the variable metric is uniformly bounded as a prior. Then we obtain the subsequential convergence of the proposed algorithm.  After employing the KL property,  the corresponding whole sequential convergence is obtained and the sequence of variable metric is proved to be bounded.  We are curious whether the boundedness of the sequence of variable metrics should be a consequence of (problem-tailored) convergence results or a prerequisite that determines the convergence,  we prefer the former.

\section{Preliminaries}
Note that all the notation is primarily taken from Rockafellar and Wets \cite{RockafellarWets2009}.  With $\mathbb R$ and $\overline{\mathbb R}:= \mathbb R \cup \{\infty\}$ we denote the real and extended real line,
respectively.  
We use $0$ to represent scalar zero, zero vector as well as zero matrix of the
appropriate dimension. 
Recall that $\mathbb R^n$ are an $n$-dimensional Euclidean space with the inner product $\left< \cdot, \cdot\right>$ and the norm denoted by $\|\cdot\|$.  
We write $A\succ 0$ ($A \succeq 0$) for $A \in \mathbb R^{n \times  n}$ if $A$ is positive (semi)definite.  We say a symmetric matrix $A$ is \emph{uniformly positive definite} if it is positive definite and there exists a positive real number $m>0$ such that $\lambda_{\min}(A)\geq m$, where $\lambda_{\min}(A)$ is the minimum eigenvalue of matrix $A$, the set of such $A$ is defined as
$$\mathbb R^{n \times n}_{\geq m}:=\left\{A\in \mathbb R^{n\times n} \,|\, \lambda_{\min}(A) \geq m, m>0\right\}.$$
We write $\|\cdot\|_A:=\sqrt{\left<A\cdot, \cdot\right>}$ for the norm induced by a given $A \succ 0$.

The effective domain of an extended real-valued function ${h}:{\R^n}\to \overline{\mathbb R}$ is denoted by $\dom h:=\{x\in\R^n \,|\, h(x) < \infty \}$.
We say that $h$ is \emph{proper} if $\dom h \neq \emptyset$
and \emph{lower semicontinuous} (lsc) if $h(\bar{x}) \leq \liminf_{x\to\bar{x}} h(x)$ for all $\bar{x} \in \R^n$.
Given a proper and lsc function ${h}: {\mathbb R} \to \overline{\mathbb R}$ and a point $\bar{x}\in\dom h$, we appeal 
to $h$-\emph{attentive} convergence of a sequence $\big(x^k \big)_{k \in \mathbb N}$:
\begin{equation}
	x^k \toattentive{h} \bar{x}
	{}\quad:\Longleftrightarrow\quad{}
	x^k \to \bar{x}
	{}\quad\text{with}\quad{}
	h(x^k) \to h(\bar{x}) .
\end{equation}
By \cite[Definition~8.3]{RockafellarWets2009}, we denote by $\hat{\partial} h: {\mathbb R^n} \to {\mathbb R^n}$ the \emph{regular subdifferential} of $h$, where
\begin{equation}
	v \in \hat{\partial} h(\bar{x})
	{}\quad:\Longleftrightarrow\quad{}
	\liminf_{\substack{x\to\bar{x}\\x\neq\bar{x}}} \frac{h(x) - h(\bar{x}) - \langle v, x-\bar{x}\rangle}{\|x-\bar{x}\|} \geq 0 .
\end{equation}
The (limiting) \emph{subdifferential} of $h$ is ${\partial h}: {\mathbb R^n} \to {\mathbb R^n}$, where $v \in \partial h(\bar{x})$ if and only if there exist sequences $\big(x^k\big)_{k \in \mathbb N}$ and $\big(v^k\big)_{k \in \mathbb N}$ such that $x^k \toattentive{h} \bar{x}$ and $v^k \in \hat{\partial} h(x^k)$ with $v^k \to v$.  A vector $v \in \mathbb R^n$ is a horizon subgradient of $h$ at $\bar x$, if there are sequences $x^k\toattentive{h} \bar{x}$,  $v^k \in \hat{\partial} h(x^k)$, one has $\lambda_k v_k \to v$ for some sequence $\lambda_k \searrow 0$.  The set of all horizon subgradients $\partial^{\infty} h(\bar x)$ is called \emph{horizon subdifferential}.  If $f$ is convex and differentiable at $\bar x$, then $\partial f(\bar x) =\{\nabla f(\bar x)\}$.
The subdifferential of $h$ at $\bar{x}$ satisfies $\partial(h+h_0)(\bar{x}) = \partial h(\bar{x}) + \nabla h_0(\bar{x})$ for any ${h_0}: {\mathbb R}\to {\overline {\mathbb R}}$ continuously differentiable around $\bar{x}$ \cite[Exercise~8.8]{RockafellarWets2009}.
We set $\hat{\partial} h(\bar{x}):=\partial h(\bar{x}):=\emptyset$ for each $\bar{x}\notin\dom h$ for completeness.  
With respect to the minimization of $h$, we say that $x^\ast \in \dom h$ is \emph{stationary} if $0 \in \partial h(x^\ast)$, which constitutes a necessary condition for the optimality of $x^\ast$ \cite[Theorem~10.1]{RockafellarWets2009}.

We next introduce the so-called local model function, before that, let us give the modified definition of the growth function suitable for the manuscript, which is essentially based on \cite{drusvyatskiy2021nonsmooth, ochs2019non}.
\begin{definition}\label{Def:GrowthFun}
(Growth function) 
An univariate function $\omega: \mathbb R_+ \to {\mathbb R}_+$ is called \textit{growth function} if it is differentiable and satisfies $\omega(0)={\omega}'_+(0)=0$, where $\omega'_+$ denotes the one sided (right) derivative of $\omega$. If, in addition $\omega'(t)>0$ for $t>0$ and equalities $\lim_{t \downarrow 0} \omega'(t)=\lim_{t \downarrow 0} \omega(t)/\omega'(t)=0$ hold, we say that $\omega$ is a \textit{proper growth function}.
\end{definition}
Note that \cite{drusvyatskiy2021nonsmooth} defined the growth function by requiring $\omega'(t)>0$ for all $t>0$ (i.e., $\omega$ is increasing on $(0, +\infty)$), however,  which has been relaxed by \cite{ochs2019non} and \Cref{Def:GrowthFun}. Through the growth function, an abstract description of a first-order oracle by the so-called model function is given in \cite{ochs2019non}, please also see \cite[Definition~5]{mukkamala2022global}.  Based on the growth function, we now give the definition of the local model function.
\begin{definition}\label{Def:ModelFunc}
(Local model function) Let $f$ be a proper lower semicontinuous function. A proper lower semicontinuous function $f_{\bar x}(x): \mathbb R^n \to \overline{\mathbb R}$ with $\dom f_{\bar x}=\dom f$ 
is called \textit{local model function} for $f$ around the \textit{model center} $\bar x\in \dom f$, if there exists a growth function $\omega_{\bar x}$ dependent on $\bar x$ such that 
\begin{equation}\label{Eq:def of model fun}
\forall x \ \text{approching to}\ \bar x: \qquad |f(x)-f_{\bar x}(x)| \leq \omega_{\bar x}(\norm{x-\bar x})
\end{equation}
holds. 
\end{definition}
In fact,  the growth function around the model center defined in \cref{Def:GrowthFun} approaches to the origin,  which might vary particularly rapidly away the model center. Therefore, \Cref{Def:ModelFunc} provides more freedom to choose a suitable model function,  simultaneously obeys the core rule: the model function needs to approximate the function well near the function center.  In other words,  we only need to bound the model error $|f_{\bar x}(x)-f(x)|$ for such $x$ close to $\bar x$,  however we do not mind characteristics of $f_{\bar x}(x)$ when $x$ is away from $\bar x$, where the value of the corresponding growth function might be large.  Therefore, we call the function defined in \cref{Def:ModelFunc} as \textit{local} model function.  

Obviously, for any model center $\bar x \in \dom f$,  one has
\begin{equation*}
f_{\bar x}=(f_{\bar x}-f)+f=:g_{\bar x}+f.
\end{equation*}
Then
\begin{equation*}
\partial f_{\bar x}(x) \subset \partial g_{\bar x}(x)+ \partial f(x) \quad \forall x\in \dom f
\end{equation*}
holds if $g_{\bar x}$ is smooth \cite[Exercise~8.8]{RockafellarWets2009} or the combination of $v_1 \in \partial^{\infty} g_{\bar x}$ and $v_2 \in \partial^{\infty} f$ with $v_1+v_2=0$ is unique and satisfies $v_1=v_2=0$ \cite[Corollary~10.9]{RockafellarWets2009}. Meanwhile, for the locally Lipschitz $f$, we know that $\partial^{\infty} f(x)=\{0\}$ \cite[Theorem~1.22]{Mordukhovich2018}. Motivated by these considerations, we establish the following first-order relationship, which serves as the optimality condition when the model function acts as the primary component in algorithm-driven subproblems.

\begin{proposition}\label{Prop: IRelation of subd}
Let $f$ be a proper lower semicontinuous function and $\bar x \in \dom f$ be arbitrarily fixed.  Moreover, denote $f_{\bar x}$ as the model function of $f$ at $\bar x$.  For any fixed $\tilde x \in \dom f$ and a constant $L>0$, one has
\begin{equation}\label{Eq:Subrel}
{\partial} f_{\bar x}(\tilde x) \subset {\partial} f(\tilde x)+L B_{\|\bar x-\tilde x\|}(0),
\end{equation}
provided that the following simultaneously hold:
\begin{align}
& g_{\bar x} \ \text{is smooth or} \ \partial^{\infty} f(\tilde x)=\{0\}, \label{H1}\\
&\partial g_{\bar x}({x}) \subset L  B_{\|{x}-\bar x\|}(0) \quad \forall x \ \text{in a neighborhood of}\ \bar x.\label{H2}
\end{align}
\end{proposition}
Note that \eqref{H2} is a not restrictive requirement, some examples in \cref{Prop:App1} and \cref{Prof:App2} are given for the general cases that $f$ is composite. 
\eqref{H2} covers \eqref{Eq:def of model fun} particularly when $x=\bar x$. However, the latter, to the best of our knowledge, does not allow to establish the desired first-order information even though the classical and easiest growth function being quadratic is employed.  When $x=\bar x$,  a specific case of \cref{Prop: IRelation of subd} (\eqref{H1} and \eqref{H2} are not needed any more) illustrates that 
\begin{equation}\label{Eq:Subsame}
\partial f_{\bar x} (\bar x) \subset  \partial f(\bar x),
\end{equation}
which has already been given in \cite[Lemma~A.1]{ochs2019model}.  

Our global convergence theory relies on the so-called Kurdyka-{\L}ojasiewicz property that 
plays a 
central role in our subsequent convergence analysis. The 
version stated here is a generalization of the classical 
Kurdyka-{\L}ojasiewicz inequality to nonsmooth functions as
introduced in \cite{AttouchBolteRedontSoubeyran2010,BolteDaniilidisLewisShiota2007}
and afterwards used in the local convergence analysis of several
nonsmooth optimization methods, cf.\ \cite{AttouchBolteSvaiter2013,BolteSabachTeboulle2014,BotCsetnek2016,BotCsetnekLaszlo2016,Ochs2018}
for a couple of examples.

\begin{definition}\label{Def:KL-property}
Let $ g: \mathbb R^n \to \R \cup \{ + \infty \} $ be proper and
lower semicontinuous. We say that $ g $ has the 
\emph{KL property} (Kurdyka-{\L}ojasiewicz property) 
at $ x^* \in \dom \partial g $ if there exist a constant 
$ \eta > 0 $, a neighborhood $ U $ of $ x^* $, and a continuous
concave function $ \varphi: [0, \eta] \to \R_+ $ with 
$$
   \varphi (0) = 0, \quad \varphi \in C^1 (0, \eta), \quad
   \text{and} \quad \varphi'(t) > 0 \quad \text{for all }
   t \in (0, \eta)
$$
such that the \emph{KL inequality}
$$
   \varphi ' \big( g(x) - g(x^*) \big) \dist \big( 0, \partial 
   g(x) \big) \geq 1
$$
holds for all $ x \in U \cap \big\{ x \in \mathbb R^n \mid 
g(x^*) < g(x) < g(x^*) + \eta \big\} $.
\end{definition}

The function $ \varphi $ is called the \emph{desingularization function}.
We note that there exist classes of functions where the 
KL property holds with the corresponding desingularization function
given by $ \varphi (t) := c t^{\theta} $ for 
$ \theta \in (0,1] $ and some constant $ c > 0 $, where the
parameter $ 1-\theta $ is called the \emph{{KL} exponent},
see \cite{BolteDaniilidisLewisShiota2007,Kurdyka1998}.
It is well known that classes of functions definable in an
o-minimal structure \cite{van1996geometric} have the KL property, which can be achieved 
for the sets or functions which are semialgebraic and globally subanalytic \cite{BolteDaniilidisLewisShiota2007}.

\section{Algorithm and Convergence Analysis}\label{Sec:Alg and Res}
This section aims to propose our model quasi-Newton methods, whose subproblems are the minimization of the regularized model function of the objective function $f$, and then demonstrates the convergence of the entire sequence of iterates in the presence that $f$ has the KL property at some accumulation point.  While, throughout, we do not make any boundedness assumption on the sequence of iterates.  For the convergence analysis, it is reasonable to assume that there exists at least one accumulation point, i.e., not every subsequence is bounded.
\subsection{Algorithm}
The overall method is stated in \Cref{Alg:ModelQuasiNewton}. 
\begin{algorithm}[H]\caption{Model Proximal Quasi-Newton Methods}
	\label{Alg:ModelQuasiNewton}
	\begin{algorithmic}[1]
		\REQUIRE $\tau > 1$, $\mu>0$, $0 < \gamma_{\min} \leq  \gamma_{\max} < \infty$, 
			$\delta \in (0,\frac{1}{2})$, $H:\mathbb R^n \to \mathbb R_{\geq \mu}^{n \times n}$ is continuous.
		\STATE Set $k := 0$. Choose $x^0 \in \dom f$ and set $H_0:=H(x^0)$.
		\WHILE{$x^k$ is not a stationary point of $f$,}
		\STATE \label{step:UnifPosD}Choose $ \gamma_k^0 \in [ \gamma_{\min}, \gamma_{\max}] $,  update the variable metric $H_k:=H(x^k)$.
		\STATE\label{step:subproblem_solve_ModelQuasiNewton} 
			For $ i = 0, 1, 2, \ldots $, compute a solution $ x^{k,i} $ of
      		\begin{equation}\label{Eq:NonSubki}
         		\min_x \ f_{x^k}(x)+ \frac{\gamma_{k,i}}{2} \| x - x^k \|_{H_k}^2 
      		\end{equation}
      		with $ \gamma_{k,i} := \tau^{i+1} \gamma_k^0 $, until the acceptance criterion
      		\begin{equation}\label{Eq:NonStepCrit}
         		|f(x^{k,i})-f_{x^k}(x^{k,i})| \leq 
         		 \delta \frac{\gamma_{k,i}}{2}\| x^{k,i} - x^k \|_{H_k}^2 
      		\end{equation}
      		holds.
		\STATE \label{item:remark} Denote $ i_k := i $ as the terminal value, and set $ \gamma_k := 
      			\gamma_{k,i_k} $ and $ x^{k+1} := x^{k,i_k} $.
                                    \STATE \label{Eq:updateH}
                  Set $ k \leftarrow k + 1 $.
		\ENDWHILE
		\RETURN $x^k$.
	\end{algorithmic}
\end{algorithm}

In the remaining parts, we assume that $\omega_{x^k}:\mathbb R_{+} \to \mathbb R_{+}$ is the growth function
such that
\begin{equation*}
|f_{x^k}(x)-f(x)| \leq \omega_{x^k}(\|x^k-x\|) \quad \forall x \ \text{approaching to } x^k
\end{equation*}
for all $k\in \mathbb N$.

In order to guarantee the convergence, we require some technical assumptions. 
\begin{assumption}\label{Ass:ModelNewton}
\leavevmode
\begin{enumerate}[(a)]
   \item \label{item: lowerbounded} $f$ is bounded from below.
   \item \label{item: mf} The model function $f_{\bar x}$ is bounded from below by an affine function 
for all $\bar x\in \dom f$. 
     \item \label{item: Unipd} $H:\mathbb R^n \to \mathbb R_{\geq \mu}^{n \times n}$ is continuous.
    \item \label{item: existence of accumulations} There exist at least one accumulation point $x^*\in \mathbb R^n$ of the iterative sequence $\big(x^k\big)_{k\in \mathbb N}$ generated by \Cref{Alg:ModelQuasiNewton}.

\end{enumerate}
\end{assumption} 
Note that \Cref{Ass:ModelNewton} \ref{item: lowerbounded} is widely used  to guarantee that \eqref{Eq:P} is solvable.  
We require that $f_{x^k}$ can be bounded from below by an affine function in \ref{item: mf} and $H_k:=H(x^k)$ is uniformly positive definite deduced by \ref{item: Unipd}, which are proposed to guarantee that subproblems \eqref{Eq:NonSubki} in Step \ref{step:subproblem_solve_ModelQuasiNewton} for fixed $k,i \in \mathbb N$,  are coercive,  and therefore always attain a solution $x^{k,i}:=x^{k+1}$, which is actually not unique.  In addition, to ensure that \Cref{Alg:ModelQuasiNewton} is well-defined, in other words, its inner loop must terninate in finite steps.   An obvious examples about \ref{item: Unipd} are that $H$ is chosen as the Hessian operator when $f$ is $\mathcal C^2$.  Also,  $H$ can be regarded as the Hessian operator of the $\mathcal C^2$ function when $f$ is composite. We will exploit more flexible structure of $H$ in the view of geometric variational analysis and/or deep learning. Due to \ref{item: existence of accumulations}, we do not need any assumptions on the boundedness of the iterative sequence. 

Step \ref{step:subproblem_solve_ModelQuasiNewton} of \Cref{Alg:ModelQuasiNewton} contains the main computational cost since we have to solve the subproblem at each iteration, which encodes that $x^{k,i}:=x^{k+1}$, as the solution of  \eqref{Eq:NonSubki}, at least reduces the value of objective function compared with $x^k$ (See \Cref{Pro:decreasing obj}).  Moreover, Step \ref{step:subproblem_solve_ModelQuasiNewton} essentially minimizes the regularized model function of the objective function,  we need to notice that the corresponding approximation error should be small, at least around the stationary point of \eqref{Eq:P}. 

Apart from the basic \Cref{Ass:ModelNewton}, we still need \eqref{H1} and \eqref{H2} to obtain the first-order information between the objective function and its corresponding models as presented in \cref{Prop: IRelation of subd}. Here, we reformulate them as
\begin{assumption}\label{Ass:H}
For all $\bar x \in \dom f$, $g_{\bar x}:=f_{\bar x}-f$:
\begin{enumerate}[(H\arabic*)]
  \item \label{Ass: H1} $g_{\bar x}$ is smooth or $f$ is locally Lipschitz continuous.
  \item \label{Ass: H2} $\partial g_{\bar x}(\cdot) \subset L B_{\|\cdot-\bar x\|}(0)$ holds with some constant $L>0$.
\end{enumerate}
\end{assumption}

We now illustrate that the stepsize rule in Step \ref{step:subproblem_solve_ModelQuasiNewton} of \Cref{Alg:ModelQuasiNewton} is always finite.
\begin{lemma}\label{Lem:well-defined inner loop}
Let $k$ be a fixed iteration of \Cref{Alg:ModelQuasiNewton}, assume that $x^k$ is not a stationary point of \eqref{Eq:P}, and suppose that \Cref{Ass:ModelNewton} holds. Then, the inner loop in Step \ref{step:subproblem_solve_ModelQuasiNewton} of \Cref{Alg:ModelQuasiNewton} terminates in a finite number of steps.
\end{lemma}
\begin{proof}
Suppose that the inner loop of  \Cref{Alg:ModelQuasiNewton} does not terminate after a finite number of steps at iteration $k$, i.e.,  $\gamma_{k,i} \to \infty$ for $i \to \infty$. Recall that $x^{k,i}$ is a solution of \eqref{Eq:NonSubki}, which implies 
\begin{equation}\label{Eq:solution of subproblems}
f_{x^k}(x^{k,i})+\frac{\gamma_{k,i}}{2}\|x^{k,i}-x^k\|^2_{H_k} \leq f(x^k).
\end{equation}
Therefore, we have $\|x^{k,i}-x^k\|_{H_k} \to 0$ for $i\to \infty$, otherwise the left-hand side of \eqref{Eq:solution of subproblems} will go to infinity and hence be unbounded by $\gamma_{k,i} \to \infty \ (i\to \infty)$, which violates the assumption that $f$ is bounded from below in view of  \Cref{Ass:ModelNewton}\ \ref{item: lowerbounded}. Furthermore, $\|x^{k,i}-x^k\| \to 0$ is valid, hence $x^{k,i} \to x^k$ as $i\to \infty$ holds from \Cref{Ass:ModelNewton}\ \ref{item: Unipd}.  Note that the model function $f_{x^k}$ is lower semicontinuous by \Cref{Def:ModelFunc},  then taking the limit $i \to \infty$ in \eqref{Eq:solution of subproblems} yields
\begin{equation*}
f(x^k)=f_{x^k}(x^k)\leq \liminf_{i \to \infty} f_{x^k}(x^{k,i})\leq \limsup_{i \to \infty} f_{x^k}(x^{k,i}) \leq f(x^k),
\end{equation*}
where the final inequality is the consequence of \eqref{Eq:solution of subproblems}. Therefore, we have
\begin{equation}\label{Eq:AttentiveModelf}
f_{x^k}(x^{k,i}) \to f_{x^k}(x^k) \quad \text{as}\ i\to \infty.
\end{equation}

We claim that 
\begin{equation}\label{Eq:forCliam}
\liminf_{i \to \infty} \gamma_{k,i}\|x^{k,i}-x^k\|_{H_k}>0.
\end{equation}
Assume, by contradiction, that there exists a subsequence $i_l \to \infty$ such that
\begin{equation}\label{Eq:gamma times difference of sequence equals to 0}
\liminf_{l \to \infty} \gamma_{k,i_l} \|x^{k,i_l}-x^k\|_{H_k}=0.
\end{equation}
Since $x^{k,i_l}$ is a solution of \eqref{Eq:NonSubki},  one has
\begin{equation}\label{Eq:optimal condition for subproblems}
0 \in {\partial} f_{x^k}(x^{k, {i_l}})+\gamma_{k, i_l}H_k(x^{k,i_l}-x^k).
\end{equation}
Taking the limit $l\to \infty$ in \eqref{Eq:optimal condition for subproblems}, combined with \eqref{Eq:AttentiveModelf} and \eqref{Eq:gamma times difference of sequence equals to 0},  implies
\begin{equation*}
0\in \partial f_{x^k}(x^k) \subset \partial f(x^k)
\end{equation*}
from \eqref{Eq:Subsame}.  Therefore,  $x^k$ is a stationary point of \eqref{Eq:P}. That is a contradiction, hence \eqref{Eq:forCliam} holds.  In view of \Cref{Ass:ModelNewton}\ \ref{item: Unipd}, one has $\|x^{k,i}-x^k\|_{H_k} \geq \sqrt{\mu}\|x^{k,i}-x^k\|$ for the fixed $k$.
Therefore, \eqref{Eq:forCliam} and the fact that $x^{k,i} \to x^k $ as $i\to \infty$ imply that there exists some $\rho_k >0$ satisying $\gamma_{k,i}\|x^{k,i}-x^k\|_{H_k} \geq \rho_k$, and hence
\begin{equation}\label{Eq: boundedness of growth function}
\delta\frac{\gamma_{k,i}}{2}\|x^{k,i}-x^k\|^2_{H_k}\geq \frac{\delta\rho_k}{2}\|x^{k,i}-x^k\|_{H_k}\geq \frac{\delta \rho_k \sqrt{\mu}}{2}\|x^{k,i}-x^k\| \geq o(\|x^{k,i}-x^k\| )
\end{equation}
holds for sufficiently large $i$ and the fixed $k$. Hence, \eqref{Eq: boundedness of growth function} yields for sufficiently large $i$,
\begin{equation*}
|f(x^{k,i})-f_{x^k}(x^{k,i})| \leq \omega_{x^k}(\|x^{k,i}-x^k\|)=o(\|x^{k,i}-x^k\|) \leq \delta\frac{\gamma_{k,i}}{2}\|x^{k,i}-x^k\|^2_{H_k},
\end{equation*}
 which contradicts $\gamma_{k,i} \to \infty$ and validates Step \ref{step:subproblem_solve_ModelQuasiNewton} of \Cref{Alg:ModelQuasiNewton} for finite $\gamma_{k,i}$.
\end{proof}

In the following, we prove that the sequence of objective values is decreasing and also convergent, which plays an central role for the convergence analysis.
\begin{proposition}\label{Pro:decreasing obj}
Let \Cref{Ass:ModelNewton} hold.  Suppose that the sequence $\big(x^k\big)_{k\in \mathbb N}$ is generated by \Cref{Alg:ModelQuasiNewton}, then $\big(f(x^k)\big)_{k\in \mathbb N}$ is a decreasing sequence and $\|x^{k+1}-x^k\| \to 0$ holds. 
\end{proposition}
\begin{proof}
Using \eqref{Eq:NonSubki} and \eqref{Eq:NonStepCrit}, we have
\begin{equation}\label{Eq:distance of decreasing for obj}
\begin{aligned}
f(x^{k+1})-f(x^k) &=f(x^{k+1})-f_{x^k}(x^{k+1})+f_{x^k}(x^{k+1})-f(x^k)\\
& \leq \delta \frac{\gamma_k}{2}\|x^{k+1}-x^k\|^2_{H_k}-\frac{\gamma_k}{2}\|x^{k+1}-x^k\|_{H_k}^2\\
&=-(1-\delta)\frac{\gamma_k}{2}\|x^{k+1}-x^k\|^2_{H_k}\leq 0
\end{aligned}
\end{equation}
for all $k\in \mathbb N$, where the last inequality is from $\delta \in (0,1)$ and the positive definiteness of $H_k$ from \Cref{Ass:ModelNewton}\ \ref{item: Unipd}.

Since the sequence $\big(f(x^k)\big)_{k\in \mathbb N}$ is monotonically decreasing, then $\big(x^k\big)_{k\in \mathbb N} \subset \mathcal L_f(x^0):=\left\{x\in \mathbb R^n \,|\, f(x) \leq f(x^0)\right\} \subset \dom f$.  Since $f$ is bounded below in view of \Cref{Ass:ModelNewton}\ \ref{item: lowerbounded}, then taking the summation $\sum_{k=0}^{\infty}$ in \eqref{Eq:distance of decreasing for obj} implies that
\begin{equation}\label{Eq:longSeqto0}
\frac{\gamma_k}{2}\|x^{k+1}-x^k\|^2_{H_k}\to 0 \quad \text{as} \ k\to \infty.
\end{equation}
Note that $\gamma_k \geq \gamma_{\min}>0$ and \Cref{Ass:ModelNewton}\ \ref{item: Unipd}, which implies
\begin{equation*}
\frac{\gamma_{\min}\mu}{2}\|x^{k+1}-x^k\|^2 \to 0\quad \text{as} \ k\to \infty,
\end{equation*}
and therefore one has $\|x^{k+1}-x^k\| \to 0$. 
\end{proof}

\begin{proposition}\label{Prop:Convergentf}
Let \Cref{Ass:ModelNewton} hold and the sequence $\big(x^k\big)_{k \in \mathbb N}$ be generated by \Cref{Alg:ModelQuasiNewton}.  One has $f(x^k) \to f^*$ with $ f^*\geq f(x^*)$ holds.
\end{proposition}
\begin{proof}
Let $\big(x^k\big)_{k\in K}$ be the subsequence convergent to $x^*$.  Furthermore $\big(x^{k+1}\big)_{k \in K}$ also converges to $x^*$ by \Cref{Pro:decreasing obj}.  Since $f$ is lower semicontinuous, we have
\begin{equation}\label{Eq:lsc}
f(x^*) \leq \liminf_{k \to_K \infty} f(x^{k+1}).
\end{equation}
On the other hand, by \Cref{Pro:decreasing obj}, the entire sequence $\big(f(x^k)\big)_{k \in \mathbb N}$ is monotonically decreasing.  Since it is also bounded from below by \Cref{Ass:ModelNewton} \ref{item: lowerbounded},  the whole sequence $\big(f(x^k)\big)_{k \in \mathbb N}$ converges and we denote its limit as $f^*$.  Obviously, $f^*\geq f(x^*)$.
\end{proof}

\subsection{Subsequential Convergence Analysis}

For the subsequential convergence analysis,  \cref{Ass:H} should be employed. 
\begin{proposition}\label{Prof:boundedSubGamma}
Let \Cref{Ass:ModelNewton} and \Cref{Ass:H} \ref{Ass: H2} 
 hold and the sequence $\big(x^k\big)_{k\in \mathbb N}$ be generated by  \Cref{Alg:ModelQuasiNewton},  and let $\big(x^k\big)_{k\in  K}$ be a subsequence converging to the point $x^*$.  Then $\gamma_k\|x^{k+1}-x^k\| \to_K 0$ holds.
\end{proposition}
\begin{proof}
If the subsequence $\left(\gamma_k\right)_{k\in K}$ is bounded, the statement holds by \cref{Pro:decreasing obj}.  It remains to consider the case where the subsequence is unbounded.  Without loss of generality, we may assume that $\gamma_k \to_K \infty$ and the acceptance criterion \eqref{Eq:NonStepCrit} is violated in the first iteration of the inner loop for each $k \in \mathbb N$. Then, for $\hat{\gamma}_k:=\gamma_k/ \tau$, we also have $\hat{\gamma}_k \to_K \infty$, the corresponding vector $\hat{x}^k:=x^{k,i_k-1}$ does not satisfy the stepsize condition from \eqref{Eq:NonStepCrit}, i.e., we have
\begin{equation}\label{Eq:ViolatedAccp}
\left|f(\hat{x}^k)-f_{x^k}(\hat x^k)\right| >\delta \frac{\hat{\gamma}_k}{2}\|\hat{x}^k-x^k\|_{H_k}^2 \quad \forall k\in K,
\end{equation}
which implies that $\hat x^k \neq x^k$ for all $k\in K$.
Meanwhile, since $\hat x^k$ solves the corresponding subproblems \eqref{Eq:NonSubki} with $\hat{\gamma}_k$, so, we have
\begin{equation}\label{Eq:SolutionofHatx}
f_{x^k}(\hat{x}^k)+\frac{\hat{\gamma}_k}{2}\|\hat{x}^k-x^k\|_{H_k}^2 \leq f(x^k)  \leq f(x^0)\quad \forall k\in K,
\end{equation}
where the second inequality is obtained because $\big(f(x^k)\big)_{k\in \mathbb N}$ is decreasing.

On the other hand,  exploiting the fact that $x^{k+1}$ and $\hat x^k$ are solutions of subproblems \eqref{Eq:NonSubki} with parameters $\gamma_k$ and $\hat \gamma_k$,  we find
\begin{equation}\label{Eq:ControlledbyModel}
\begin{aligned}
f_{x^k}(x^{k+1})+\frac{\gamma_k}{2}\|x^{k+1}-x^k\|_{H_k}^2 &\leq f_{x^k}(\hat x^k)+\frac{\gamma_k}{2}\|\hat x^{k}-x^k\|_{H_k}^2,  \\
f_{x^k}(\hat x^k)+\frac{\hat \gamma_k}{2}\|\hat x^{k}-x^k\|_{H_k}^2 &\leq f_{x^k}(x^{k+1})+\frac{\hat \gamma_k}{2}\|x^{k+1}-x^k\|_{H_k}^2,
\end{aligned}
\end{equation}
for all $k \in K$.  Adding these two inequalities and noting that $\gamma_k:=\tau \hat \gamma_k$ imply that $\|x^{k+1}-x^k\|^2_{H_k} \leq \|\hat x^k-x^k\|^2_{H_k}$ for all $k \in K$.  Therefore,  we obtain from the second inequality in \eqref{Eq:ControlledbyModel} that
\begin{equation}\label{Eq:ControlledbyModel2}
f_{x^k}(\hat x^k) \leq f_{x^k}(x^{k+1}) \quad \forall k\in K.
\end{equation}
By \Cref{Pro:decreasing obj},  \Cref{Def:GrowthFun}, and \eqref{Eq:eta},  we have
\begin{equation}\label{Eq:Model}
f_{x^k}(x^{k+1}) \leq f(x^{k+1})+\omega_{x^k} (\|x^{k+1}-x^k\|) \leq f(x^0)+\omega_{x^k}(\|x^{k+1}-x^k\|)<\infty 
\end{equation}
for all $k \geq \hat k$ and $k \in K$, in other words, $f_{x^k}(x^{k+1})$ is finite for all $k \geq \hat k$ and $k\in K$.
Therefore,  \eqref{Eq:ControlledbyModel2} implies that $f_{x^k}(\hat x^k)$ is finite for all $\hat k \leq k \in K$.  Hence, let us look at \eqref{Eq:SolutionofHatx} again, we can definitely say that $\|\hat x^k-x^k\|_{H_k} \to_K 0$ by $\hat \gamma_k \to_K \infty$ (otherwise the left-hand side of \eqref{Eq:SolutionofHatx} goes to infinity),  which implies that $\|\hat x^k -x^k\| \to_K 0$ by \Cref{Ass:ModelNewton}\ \ref{item: Unipd}.

From \eqref{Eq:ViolatedAccp} and \Cref{Ass:ModelNewton}\ \ref{item: Unipd}, we obtain
\begin{equation}\label{Eq:ContrforGamma}
 \delta \frac{\hat{\gamma}_k}{2}\|\hat{x}^k-x^k\|_{H_k}^2 < |f(\hat x^k)-f_{x^k}(\hat x^k)| =|g_{x^k}(\hat x^k)|\quad \forall k\in K.
\end{equation}
On the other hand,  in view of \Cref{Ass:H} \ref{Ass: H2}, we have $\hat \partial g_{x^k}(x^k)=\{0\}$ and $\|\eta^k\| \leq L\|\hat x^k-x^k\| \ \forall \eta^k\in \hat \partial g_{x^k}(\hat x^k)$ hold for any fixed $k\in \mathbb N$.  By the definition of the regular subdifferential, we have
\begin{equation*}
\begin{aligned}
&\liminf_{\substack{x\to \hat x^k\\x\neq {\hat x^k}}} \frac{g_{x^k}(x) - g_{x^k}(\hat x^k) - \langle\eta^k, x-\hat x^k\rangle}{\|x-\hat x^k\|} \geq 0,\\
&\liminf_{\substack{x\to x^k\\x\neq x^k}} \frac{g_{x^k}(x) - g_{x^k}(x^k)}{\|x-x^k\|} \geq 0 .
\end{aligned}
\end{equation*}
Due to $\|x^k-\hat x^k\| \to_K 0$ and $x^k \neq \hat x^k$ for any $k\in K$,  then for arbitrary $\varepsilon>0$, and the large enough $k \in K$, we have 
\begin{equation*}
\begin{aligned}
&g_{x^k}(x^k) - g_{x^k}(\hat x^k) - \langle\eta^k, x^k-\hat x^k\rangle \geq -\varepsilon\|x^k-\hat x^k\|,\\
& g_{x^k}(\hat x^k) - g_{x^k}(x^k) \geq -\varepsilon\|x^k-\hat x^k\|.
\end{aligned}
\end{equation*}
They imply that
\begin{equation*}
-\langle\eta^k, x^k-\hat x^k\rangle+\varepsilon\|x^k-\hat x^k\| \geq g_{x^k}(\hat x^k) - g_{x^k}(x^k) \geq -\varepsilon \|x^k-\hat x^k\|
\end{equation*}
for those $k\in K$ large enough.  Recall again $\|\eta^k\| \leq L\|\hat x^k-x^k\|$,  it implies that
\begin{equation*}
|g_{x^k}(\hat x^k)|= | g_{x^k}(\hat x^k) - g_{x^k}(x^k)| \leq L \|\hat x^k-x^k\|^2+\varepsilon\|\hat x^k-x^k\|
\end{equation*}
holds for the sufficiently large $k \in K$.  It, together with \eqref{Eq:ContrforGamma}, implies that
\begin{equation*}
\begin{aligned}
\delta \frac{\hat{\gamma}_k}{2}\|\hat{x}^k-x^k\|_{H_k}^2 &< L \|\hat x^k-x^k\|^2+\varepsilon\|\hat x^k-x^k\|\\
& =\frac{1}{\min\{\mu, \sqrt{\mu}\}}\min\{\mu, \sqrt{\mu}\}\big(L \|\hat x^k-x^k\|^2+\varepsilon\|\hat x^k-x^k\|\big)\\
& \leq \frac{1}{\min\{\mu, \sqrt{\mu}\}} \big(L\mu\|\hat x^k-x^k\|^2+\varepsilon \sqrt{\mu}\|\hat x^k-x^k\|\big)\\
& \leq \frac{1}{\min\{\mu, \sqrt{\mu}\}} \big(L\|\hat x^k-x^k\|_{H_k}^2+\varepsilon \|\hat x^k-x^k\|_{H_k}\big)
\end{aligned}
\end{equation*}
for the sufficiently large $k\in K$. Recall again that $\hat x^k \neq x^k$ for all $k\in K$, we have 
\begin{equation*}
\delta \frac{\hat{\gamma}_k}{2}\|\hat{x}^k-x^k\|_{H_k}<\frac{1}{\min\{\mu, \sqrt{\mu}\}}\big(L \|\hat x^k-x^k\|_{H_k}+\varepsilon\big)
\end{equation*}
for all sufficiently large $k\in K$. Recall again that 
$\|\hat x^k-x^k\|_{H_k}\to_{K} 0$, it implies that $\hat \gamma_k\|\hat{x}^k-x^k\|_{H_k} \to_{ K}0$.
By $\|x^{k+1}-x^k\|_{H_k} \leq \|\hat x^k-x^k\|_{H_k}$, we have
\begin{equation*}
\gamma_k\|x^{k+1}-x^k\|\leq \gamma_k\frac{1}{\sqrt{\mu}}\|x^{k+1}-x^k\|_{H_k}\leq \frac{\tau }{\sqrt{\mu}}\hat \gamma_k\|\hat x^{k}-x^k\|_{H_k}\to_{ K}  0. 
\end{equation*}
This completes the proof.
\end{proof}

The following is the main (subsequential) convergence result of \Cref{Alg:ModelQuasiNewton}.
\begin{theorem}\label{Th:SubConv}
Let \Cref{Ass:ModelNewton} and \Cref{Ass:H} \ref{Ass: H2}
hold, $f$ be further locally Lipschitz continuous,  the sequence $\big(x^k\big)_{k\in \mathbb N}$ be generated by  \Cref{Alg:ModelQuasiNewton},  and let $\big(x^k\big)_{k\in  K}$ be a subsequence converging to the point $x^*$.  Then, $x^*$ is a stationary point of \eqref{Eq:P}.
\end{theorem}
\begin{proof}
Since $x^{k+1}$ is a solution of subproblems \eqref{Eq:NonSubki}, then one has
\begin{equation}\label{Eq:forSubfx}
0\in {\partial} f_{x^k}(x^{k+1})+\gamma_kH_k(x^{k+1}-x^k) \quad \forall k \in \mathbb N.
\end{equation}
We know \Cref{Ass:H} holds since $f$ is locally Lipschitz continuous and \Cref{Ass:H} \ref{Ass: H2} is required, recall again \cref{Prop: IRelation of subd}, one has 
\begin{equation}\label{Eq:forConv}
\gamma_kH_k(x^{k}-x^{k+1}) \in {\partial} f_{x^k}(x^{k+1}) \subset \partial f(x^{k+1})+L B_{\|x^{k+1}-x^{k}\|}(0),
\end{equation}
 for all $k\in \mathbb N$.
 Hence,  by \Cref{Prop:BoundedGamma}, we have
\begin{equation*}
\begin{aligned}
\dist(0,\partial f(x^{k+1}))&  \leq {\gamma_k}\|H_k(x^{k}-x^{k+1})\|+L\|x^{k+1}-x^k\|\\
& \leq \gamma_k\|H_k\|\|x^{k}-x^{k+1}\|+L\|x^{k+1}-x^k\|
\end{aligned}
\end{equation*}
for all $k \in \mathbb N$.  By \Cref{Ass:ModelNewton}~\ref{item: Unipd} as well as the fact that $x^k \to_K x^*$ and $\|x^{k+1}-x^k\| \to 0$ in the view of \cref{Pro:decreasing obj}, we have that $\|H_k\| \to_K \|H(x^*)\|$ and $\gamma_k\|x^{k+1}-x^k\| \to_K 0$ deduced by \cref{Prof:boundedSubGamma}. Meanwhile, one has $f(x^{k+1}) \to_K f(x^*)$ by the local Lipschitz continuity of $f$.  Therefore $0 \in \partial f(x^*)$ holds, i.e., $x^*$ is a stationary point of \eqref{Eq:P}.
\end{proof}
\subsection{Sequential Convergence Analysis}
\cref{Th:SubConv} illustrates that any cluster point of the sequence generated by \Cref{Alg:ModelQuasiNewton} is stationary. In order to obtain the corresponding convergence result of the whole generated sequence,  we need to assume the objective function $f$ has the Kurdyka-{\L}ojasiewicz property at the accumulation point.
\begin{assumption}\label{Ass:KL}
\leavevmode
\begin{enumerate}[(a)]
 \item \label{item: KL} $f$ has the KL property at $x^*$ which is from \cref{Ass:ModelNewton}  \ref{item: existence of accumulations}.
\end{enumerate}
\end{assumption}
By employing \Cref{Ass:ModelNewton}, \Cref{Ass:H}, and \Cref{Ass:KL},  we first show that the stepsize is bounded whenever the iterates stay in some neighborhood centered around the acccumulation point $x^*$. The result is used to guarantee the \textit{local} relative error condition holds on this neighborhood in \Cref{Lem:BoundedSubf}, which is mostly necessary as a proof technique when using the \rm{KL} property.  When the KL property of $f$ is assumed in \Cref{Thm:ConvofEntireSeq}, we prove that the iterates with sufficiently large counts, in turn,  stay in this neighborhood (therefore, the corresponding stepsize is bounded),  and that the sequence generated by \Cref{Alg:ModelQuasiNewton} has a finite length and is consequently convergent to the stationary point.

Before declaring these results, let us give some notation for the convenience.  Let sufficiently small (see \cref{Pro:decreasing obj}) $\eta>0$ be the corresponding constant of the associated desingularization function $\varphi$ in \Cref{Def:KL-property} and $\hat k \in \mathbb N$ be a sufficiently large index such that
\begin{equation}\label{Eq:eta}
\sup_{k \geq \hat k} \|x^{k+1}-x^k\| \leq \eta.
\end{equation}
We set $\rho:=\eta+\frac{1}{2}$
and define
the index set
\begin{equation}\label{Eq:ind}
I_{\rho}:=\left\{k\in \mathbb N \,|\, x^k \in B_{\rho}(x^*)\right \},
\end{equation}
as well as the compact set
\begin{equation}\label{Eq:Cset}
C_{\rho}:=B_{\rho}(x^*) \cap \mathcal L_{f}(x^0),
\end{equation}
where $\mathcal L_{f}(x^0):=\left\{x\in \mathbb R^n \,|\, f(x) \leq f(x^0)\right\}$ is the sublevel set of $f$ with respect to $x^0$, the starting point exploited in \Cref{Alg:ModelQuasiNewton}. 

Based on these notation, we first illustrate that the sequence of stepsize is (uniformly) bounded on any bounded set.
\begin{proposition}\label{Prop:BoundedGamma}
Let \Cref{Ass:ModelNewton} and \Cref{Ass:H} \ref{Ass: H2} 
 hold, $g_{\bar x}$ is smooth for all $\bar x \in \dom f$, and the sequence $\big(x^k\big)_{k\in \mathbb N}$ be generated by  \Cref{Alg:ModelQuasiNewton}, then for $\rho$ defined after \eqref{Eq:eta},  there exists some constant $\bar \gamma_{\rho}>0$ dependent on $\rho$ such that $\gamma_k \leq \bar \gamma_{\rho}$ holds for all $k\in I_{\rho} $, where $I_{\rho}$ is denoted in \eqref{Eq:ind}.
\end{proposition}
\begin{proof}
Now assume, by contradiction, that there is a subsequence $\big(\gamma_k\big)_{k\in K}$ with $x^k \in B_{\rho}(x^*)$ for all $k \in K$ such that $\big(\gamma_k\big)_{k\in K}$ is unbounded.  Without loss of generality, we may assume that $\gamma_k \to_K \infty$ and the acceptance criterion \eqref{Eq:NonStepCrit} is violated in the first iteration of the inner loop for each $k \in \mathbb N$. 

Let us denote $\hat{\gamma}_k :=\gamma_k/\tau$, the corresponding vector as $\hat{x}^k:=x^{k,i_k-1}$. By the proof in \cref{Prof:boundedSubGamma}, we obtain
\begin{equation*}
\delta \frac{\hat{\gamma}_k}{2}\|\hat{x}^k-x^k\|_{H_k}^2<|g_{x^k}(\hat x^k)|=|g_{x^k}(\hat x^k)-g_{x^k}(x^k)|
\end{equation*}
for all $k\in K$. Recall that \Cref{Ass:H} \ref{Ass: H2} and $g_{x^k}$ is smooth,  by the mean-value theorem, then there exists a vector $\xi^k\in \mathbb R^{n}$ on the segment between $x^k$ and $\hat x^k$ satisfying
\begin{equation*}
\begin{aligned}
\delta \mu\frac{\hat{\gamma}_k}{2}\|\hat{x}^k-x^k\|^2& \leq \delta \frac{\hat{\gamma}_k}{2}\|\hat{x}^k-x^k\|_{H_k}^2< \big|\left<\nabla g_{x^k}(\xi^k), \hat x^k-x^k\right>\big| \\
&\leq \|\nabla g_{x^k}(\xi^k)\|\|\hat x^k-x^k\|
\leq L\|\xi^k-x^k\|\|\hat x^k-x^k\|\\
& \leq L\|\hat x^k-x^k\|^2 \quad \forall k\in K.
\end{aligned}
\end{equation*}
Recall again the fact that $\hat x^k\neq x^k$ for all $k \in K$, then the subsequence $\big(\hat \gamma_k\big)_{k\in K}$ must be bounded, which,  in turn, implies the boundedness of the subsequence $\big(\gamma_k\big)_{k\in K}$, contradicting our assumption. This completes the proof.
\end{proof}
Note that the requirement in \Cref{Ass:H} \ref{Ass: H2} can be relaxed as follows: Let $L_{\bar x}$ be the constant dependent on $\bar x$ satisfying \ref{Ass: H2}.  Accordingly,  from the respective of \Cref{Alg:ModelQuasiNewton}, define $L:=\sup_{x^k\in B_{\eta}(x^*)} L_{x^k}$ in \cref{Prop:BoundedGamma} and throughout the subsequent proof. 

We next give the following weaker version of the classical relative error condition,  which is necessay when the KL property is used.  We omit the corresponding proof since it is highly similar with \cref{Th:SubConv}.
\begin{lemma}\label{Lem:BoundedSubf}
Let \Cref{Ass:ModelNewton}, \Cref{Ass:H} with smooth $g_{\bar x}$ for all $\bar x \in \dom f$, and \Cref{Ass:KL} hold,  the sequence $\big(x^k\big)_{k \in \mathbb N}$ be generated by \Cref{Alg:ModelQuasiNewton}.  Then there exists a constant $L>0$ such that 
\begin{equation*}
\dist\big(0, \partial f (x^{k+1})\big) \leq  \bar \gamma_{\rho}\|H_k\|\|x^{k+1}-x^k\|+L\|x^{k+1}-x^k\|
\end{equation*}
holds for all sufficiently large $k \geq \hat k$ and $k \in I_{\rho}$, where $\bar \gamma_{\rho}, I_{\rho}$ are denoted in \Cref{Prop:BoundedGamma} and \eqref{Eq:ind}, respectively.
\end{lemma}
Based on the weaker relative error condition, we next illustrate that the whole sequence generated by \Cref{Alg:ModelQuasiNewton} is convergent to a stationary point.
\begin{theorem}\label{Thm:ConvofEntireSeq}
Let \Cref{Ass:ModelNewton}, \Cref{Ass:H} with smooth $g_{\bar x}$ for all $\bar x \in \dom f$, and \Cref{Ass:KL} hold,  the sequence $\big(x^k\big)_{k \in \mathbb N}$ be generated by \Cref{Alg:ModelQuasiNewton}. Then $\big(x^k\big)_{k\in \mathbb N}$ converges to $x^*$, and $x^*$ is a stationary point.
\end{theorem}
\begin{proof}
We know that the whole sequence $\big(f(x^k)\big)_{k\in \mathbb N}$ is monotonically decreasing and convergent to $f^*\geq f(x^*)$ in view of \cref{Prop:Convergentf}. This implies that $f(x^k) \geq f(x^*)$ holds for all $k\in \mathbb N$.

Now, suppose we have $f(x^k)=f(x^*)$ for some index $k\in \mathbb N$. Then, by monotonicity, we also get $f(x^{k+1})=f(x^*)$. Consequently, we obtain from \eqref{Eq:distance of decreasing for obj} that
\begin{equation}
0\leq (1-\delta)\frac{\gamma_{k-1}}{2}\|x^{k}-x^{k-1}\|_{H_{k-1}} \leq f(x^k)-f(x^{k+1})=0.
\end{equation}
By \Cref{Ass:ModelNewton}\ \ref{item: Unipd} and $\gamma_k \geq \gamma_{\min}$, thus one has $x^{k}=x^{k-1}$. By \Cref{Ass:ModelNewton} \ref{item: existence of accumulations},  $x^*$ is an accumulation point of $\big(x^k\big)_{k\in \mathbb N}$, this implies that $x^k=x^*$ and consequently $f(x^k)=f(x^*)$ holds for all $k\in \mathbb N$ sufficiently large. In particular, we have the convergence of the entire sequence $\big(x^k\big)_{k\in \mathbb N}$ (eventually constant) to $x^*$ and $f(x^k) \to f(x^*)$ in this situation.\\
\indent For the remainder of this proof. We therefore assume that $f(x^k)>f(x^*)$ holds for all $k \in \mathbb N$.  Let $\eta>0$ be the corresponding constant from the definition of the associated desingularization function $\varphi$,  $\big(x^k\big)_{k\in K}$ be the subsequence convergent to $x^*$ and $k_0 \in K$ be the sufficiently large iteration index,  one has
\begin{equation}\label{Eq:k0forKL}
0< f(x^k)-f(x^*)\leq f(x^{k_0})-f(x^*) < \eta \quad \forall k \geq k_0.
\end{equation}
Without loss of generality, we may also assume that $k_0 \geq \hat k$ (the latter being the index defined in \eqref{Eq:eta}) and that $k_0$ is sufficiently large to satisfy
\begin{equation*}
f(x^{k_0})< f(x^*)+\eta.
\end{equation*}
Let $\varphi: [0, \eta] \to [0, \infty)$ be the desingularization function which comes along with the validity of the KL property.  Due to $\varphi(0)=0$ and $\varphi'(t)>0$ for all $t\in (0, \eta)$, we obtain 
\begin{equation*}
\varphi\big(f(x^k)-f(x^*)\big) \geq 0 \quad \forall k \geq k_0.
\end{equation*}
Meanwhile, from the continuity of $H$ as required in \Cref{Ass:ModelNewton} \ref{item: Unipd}, there exists a constant $M>0$ such that 
\begin{equation}\label{Eq:H}
\|H(x)\| \leq M \quad \forall x\in C_{\rho}.
\end{equation}
For $\hat k \leq k_0 \in \mathbb N$,  we set 
\begin{equation}\label{Eq:alpha}
\alpha:=\|x^{k_0}-x^*\|+\sqrt{\frac{8\big(f(x^{k_0})-f(x^*)\big)}{\mu(1-\delta)\gamma_{\min}}}+\frac{2\bar \gamma_{\rho}M}{(1-\delta)\mu\gamma_{\min}}\varphi\big(f(x^{k_0})-f(x^*)\big),
\end{equation} 
one has $\alpha$ is sufficiently small, and hence $\alpha<\rho$.
We now claim that the following statements hold for all $k \geq k_0$:
\begin{enumerate}[(a)]
\item \label{Itm:a}$x^k \in B_{\alpha}(x^*)$,
\item \label{Itm:b} $\|x^{k_0}-x^*\| +\sum_{i=k_0}^k \|x^{i+1}-x^i\| \leq \alpha$, which is equivalent to
\begin{equation}\label{Eq:forconv}
\sum_{i=k_0}^{k} \|x^{i+1}-x^i\| \leq \sqrt{\frac{8\big(f(x^{k_0})-f(x^*)\big)}{\mu(1-\delta)\gamma_{\min}}}+\frac{2\bar \gamma_{\rho}M}{(1-\delta)\mu\gamma_{\min}}\varphi\big(f(x^{k_0})-f(x^*)\big).
\end{equation}
\end{enumerate}
We verify these two statements jointly by induction.  For $k=k_0$, statement \ref{Itm:a} holds from the definition of $\alpha$ in \eqref{Eq:alpha}.  Furthermore,  \eqref{Eq:distance of decreasing for obj} together with the monotonicity of $\big(f(x^k)\big)_{k\in \mathbb N}$ implies
\begin{equation*}
 \sqrt{\mu}\|x^{k_0+1}-x^{k_0}\| \leq \|x^{k_0+1}-x^{k_0}\|_{H_{k_0}} \leq \sqrt{\frac{2\big(f(x^{k_0})-f(x^{k_0+1})\big)}{(1-\delta)\gamma_{\min}}} \leq \sqrt{\frac{2\big(f(x^{k_0})-f^*\big)}{(1-\delta)\gamma_{\min}}}.
\end{equation*}
In particular, this shows \ref{Itm:b} holds for $k=k_0$. Suppose that the two statements hold for some $k \geq k_0$.  Using the triangle inequality, the induction hypothesis, and the definition of $\alpha$ in \eqref{Eq:alpha}, we obtain
\begin{equation*}
\begin{aligned}
\|x^{k+1}-x^*\| &\leq \sum_{i=k_0}^k \|x^{i+1}-x^i\| +\|x^{k_0}-x^*\|\\
&\leq \sqrt{\frac{8\big(f(x^{k_0})-f(x^*)\big)}{\mu(1-\delta)\gamma_{\min}}}+\frac{2(\bar \gamma_{\rho}M+L)}{(1-\delta)\mu\gamma_{\min}} \varphi\big(f(x^{k_0})-f(x^*)\big)\\
&~~~~+\|x^{k_0}-x^*\| \\
&=\alpha, 
\end{aligned}
\end{equation*}
i.e., statement \ref{Itm:a} holds for $k+1$ in place of $k$. The verification of the induction step for \ref{Itm:b} is more involved.

To this end, 
note that \eqref{Eq:k0forKL} implies that
\begin{equation}\label{Eq:fieq}
f(x^*)<f(x^i)<f(x^*)+\eta \quad \forall i \geq k_0.
\end{equation}
Recall again \Cref{Ass:KL} as well as $x^i \in B_{\alpha}(x^*) \subset B_{\rho}(x^*)$ for all $i \in \{ k_0, k_0+1, \ldots, k\}$ by our hypothesis, hence \cref{Lem:BoundedSubf} holds and \eqref{Eq:H} is available with $x:=x^i$ for those $i$,  indicating that (after a simple index shift)
\begin{equation}\label{Eq:resultabove}
\dist \big(0, \partial f(x^{i})\big) \leq (\bar \gamma_{\rho}M+L)\|x^i-x^{i-1}\| \quad \forall i\in \{k_0+1, \ldots, k+1\}.
\end{equation}
We also have
\begin{equation}\label{Eq:KLine}
\varphi'\big(f(x^i)-f(x^*)\big) \dist \big(0, \partial f(x^{i})\big) \geq 1 \quad \forall i \in \{ k_0+1, \ldots, k+1\}.
\end{equation}
Putting \eqref{Eq:resultabove} and \eqref{Eq:KLine} together implies
\begin{equation}\label{Eq:Difchi}
\varphi'\big(f(x^i)-f^*\big) \geq \frac{1}{(\bar \gamma_{\rho}M+L)\|x^i-x^{i-1}\|}  \quad \forall i\in \{ k_0+1, k_0+2, \ldots, k+1\}.
\end{equation}
To simplify the subsequent formulas, we introduce the short hand notation
\begin{equation*}
\Delta_{i,j}:=\varphi\big(f(x^i)-f(x^*)\big)-\varphi\big(f(x^j)-f(x^*)\big)
\end{equation*}
for $i,j \in \mathbb N$. The assumed concavity of $\varphi$ then implies
\begin{equation}\label{Eq:Concachi}
\Delta_{i , i+1} \geq \varphi'\big(f(x^i)-f(x^*)\big)\big(f(x^i)-f(x^{i+1})\big).
\end{equation}
Using \eqref{Eq:distance of decreasing for obj}, \eqref{Eq:Difchi}, and \eqref{Eq:Concachi}, we therefore get
\begin{equation*}
\begin{aligned}
\Delta_{i, i+1}&\geq \varphi'\big(f(x^i)-f(x^*)\big)\big(f(x^i)-f(x^{i+1})\big) \geq \frac{f(x^i)-f(x^{i+1})}{( \bar \gamma_{\rho}M+L)\|x^i-x^{i-1}\|} \\
&\geq (1-\delta)\frac{\gamma_{\min}\|x^{i+1}-x^{i}\|^2_{H_{i}}}{2 (\bar \gamma_{\rho}M+L)\|x^i-x^{i-1}\|}\geq \beta \frac{\|x^{i+1}-x^{i}\|^2}{\|x^i-x^{i-1}\|}
\end{aligned}
\end{equation*}
for all $k_0+1\leq i \leq k+1$, where we used the constant $\beta:=\frac{(1-\delta)\mu\gamma_{\min}}{2(\bar \gamma_{\rho}M+L)}$.  
Note that $a+b \geq 2\sqrt{ab}$ holds for all real numbers $a,b\geq 0$, we therefore obtain
\begin{equation*}
\frac{1}{\beta}\Delta_{i, i+1}+\|H_{i-1}\|\|x^i-x^{i-1}\| \geq 2 \sqrt{\frac{1}{\beta} \Delta_{i, i+1}\|H_{i-1}\|\|x^i-x^{i-1}\|} \geq 2\|x^{i+1}-x^i\|
\end{equation*}
for all $i \in \{ k_0+1, k_0+2, \ldots, k+1\}$. 
Summation yields from the positive definiteness of $H_i$ that
\begin{equation*}
\begin{aligned}
2\sum_{i=k_0+1}^{k+1}\|x^{i+1}-x^i\| &\leq \sum_{i=k_0+1}^{k+1} \frac{1}{\beta}\Delta_{i, i+1}+\sum_{i=k_0+1}^{k+1}\|x^i-x^{i-1}\|\\
&=\frac{1}{\beta}\Delta_{k_0+1,k+2}+\|x^{k_0+1}-x^{k_0}\|+\sum_{i=k_0+1}^{k}\|x^{i+1}-x^i\|\\
& \leq \frac{1}{\beta}\Delta_{k_0+1,k+2}+\|x^{k_0+1}-x^{k_0}\|+\sum_{i=k_0+1}^{k+1}\|x^{i+1}-x^i\|.
\end{aligned}
\end{equation*}
Subtracting the first summand  from the right-hand side, \eqref{Eq:distance of decreasing for obj},  and using the nonnegativity as well as monotonicity of the desingularization function $\varphi$, we obtain
\begin{equation*}
\sum_{i=k_0+1}^{k+1}\|x^{i+1}-x^i\| \leq \sqrt{\frac{2\big(f(x^{k_0})-f(x^*)\big)}{(1-\delta)\gamma_{\min}}}+\frac{1}{\beta}\varphi\big(f(x^{k_0})-f(x^*)\big).
\end{equation*}
Adding the term $\|x^{k_0+1}-x^{k_0}\|>0$ to both sides and using \eqref{Eq:distance of decreasing for obj} again, we obtain
\begin{equation}\label{Eq:FinitTr}
 \sum_{i=k_0}^{k+1}\|x^{i+1}-x^i\| \leq \sqrt{\frac{8\big(f(x^{k_0})-f(x^*)\big)}{(1-\delta)\gamma_{\min}}}+\frac{1}{\beta}\varphi\big(f(x^{k_0})-f(x^*)\big),
\end{equation}
by \Cref{Ass:ModelNewton}\ \ref{item: Unipd},  \eqref{Eq:FinitTr} yields that
\begin{equation*}
\sum_{i=k_0}^{k+1}\|x^{i+1}-x^i\| \leq \frac{1}{\sqrt{\mu}}\left(\sqrt{\frac{8\big(f(x^{k_0})-f(x^*)\big)}{(1-\delta)\gamma_{\min}}}+\frac{1}{\beta}\varphi\big(f(x^{k_0})-f(x^*)\big)\right).
\end{equation*}
Hence, statement \ref{Itm:b} holds for $k+1$ in place of $k$, and this completes the induction.

This easily shows that the sequence $\big(x^k\big)_{k \in \mathbb N}$ has finite length, that is
\begin{equation*}
\sum_{k=1}^{\infty}\|x^{i+1}-x^i\| < \infty,
\end{equation*}
which indicates that $\big(x^k\big)_{k \in \mathbb N}$ is a Cauchy sequence, and hence convergent. Since we already know that $x^*$ is an accumulation point of $\big(x^k\big)_{k \in \mathbb N}$, then the entire sequence $\big(x^k\big)_{k \in \mathbb N}$ converges to $x^*$. 

We now prove $f(x^k) \to f(x^*)$ in the situation where $f(x^k) > f^*$ for all $k\in \mathbb N$.  Note that $f(x^k) \to f^* \geq f(x^*)$ by \Cref{Prop:Convergentf}, now we assume that $f^*>f(x^*)$.
From $x^k \to x^*$, one has $x^k \in B_{\rho}(x^*)$ holds for all sufficiently large $k \in \mathbb N$,  there exists some $w^{k+1} \in \partial f(x^{k+1})$ ($k$ is sufficiently large) satisfying
\begin{equation*}
\|w^{k+1}\| \leq \bar \gamma_{\rho} \|H_k(x^{k+1}-x^k)\|+L\|x^{k+1}-x^k\| \leq \big(\bar \gamma_{\rho}M+L\big)\|x^{k+1}-x^k\|,
\end{equation*}
which implies that $w^{k+1} \to 0$ from \eqref{Eq:longSeqto0}.
For such $w^{k+1}$, 
from monotone decrease of $\varphi'$, \eqref{Eq:fieq}, and \Cref{Ass:KL}, one has
\begin{equation*}
\varphi'\big(f(x^k)-f(x^*)\big)\|w^{k+1}\| \geq \varphi'\big(f(x^{k+1})-f(x^*)\big)\|w^{k+1}\| \geq 1 
\end{equation*}
for all sufficiently large $ k\in \mathbb N$,
which yields a contradiction for the sufficiently large $k \in \mathbb N$.  Hence,  $f^*=f(x^*)$, in other words, $f(x^k) \to f(x^*)$ holds. 

Recall that $x^k \to x^*$,  $\|x^{k+1}-x^k\| \to 0$, and $f(x^k) \to f(x^*)$,  then taking a limit $k \to \infty$ into \eqref{Eq:forConv} yields that $0 \in \limsup_{k \to \infty} \partial f(x^{k+1}) \subset \partial f(x^*)$, which means that
$x^*$ is a stationary point of $f$.
\end{proof}
We have obtained the sequential convergence of \Cref{Alg:ModelQuasiNewton}, i.e.,  the whole sequence generated by \Cref{Alg:ModelQuasiNewton} is convergent to a stationary point,  provided that the objective function has the KL property at the accumulation point, we next give the following convergence rate under general cases of the so-called disingularization function.  For the latest result on the superlinear convergence rate, please refer to \cite{bento2025convergence,yagishita2025proximal}
\begin{theorem}\label{Thm:Rate-of-Conv}
Let \Cref{Ass:ModelNewton}, \Cref{Ass:H} with smooth $g_{\bar x}$ for all $\bar x \in \dom f$,  and \Cref{Ass:KL} hold,  the sequence $\big(x^k\big)_{k \in \mathbb N}$ be generated by \Cref{Alg:ModelQuasiNewton}.  Then the entire sequence $\big(x^k\big)_{k \in \mathbb N}$ converges
to $ x^* $, and if the corresponding desingularization function
has the form $ \varphi (t) = c t^{1-\theta} $ for some $ c > 0 $ and $\theta \in [0,1)$, 
the following statements hold:
\begin{enumerate}[(i)]
         \item \label{conv rate equal 0} if $\theta=0$, then the sequences $\big(f(x^k)\big)_{k\in \mathbb N}$ and $\big(x^k\big)_{k \in \mathbb N}$ converge in a finite number of steps to $f(x^*)$ and $x^*$, respectively.
	\item \label{conv rate less 0.5} if  $\theta\in (0, \frac{1}{2})$, the sequences $\big(f(x^k)\big)_{k\in \mathbb N}$ and $\big(x^k\big)_{k \in \mathbb N}$ either converge in a finite number of steps, or converge superlinearly to $f(x^*)$ and $x^*$, respectively.
          \item if $\theta=\frac{1}{2}$, then
	   the sequence $\big(f(x^k)\big)_{k\in \mathbb N}$ converges 
	   $ Q $-linearly to $ f(x^*) $, and
	   the sequence  $\big(x^k\big)_{k \in \mathbb N}$ converges R-linearly to $x^* $.
         \item \label{conv rate more 0.5} if $\theta\in (\frac{1}{2},1)$, then there exist some
           positive constants $\eta_1$ and $\eta_2$ such that
           \begin{align*}
           f(x^k)-f(x^*) &\leq \eta_1 k^{-\frac{1}{1-2\theta}},\\
          \|x^k-x^*\| &\leq \eta_2 k^{-\frac{\theta}{1-2\theta}}
             \end{align*}
           for sufficiently large $k$.
         
\end{enumerate}
\end{theorem}

\section{Examples}\label{Sec:Ex}
In this section,  we consider some instances of \eqref{Eq:P}, in particular, (additive) composite problems, in order to illustrate that our problem setting is much more general and hence \Cref{Alg:ModelQuasiNewton} corresponds to classical method by defining suitable model functions.

\subsection{Additive composite problems}
We consider the following (nonconvex) additive composite problem:
\begin{equation}\label{Eq:ACP}
	\min_x \ f(x) :=q(x)+h(x),
\end{equation}
where $h:\mathbb R^n \to \mathbb R$ is continuously differentiable and $q:\mathbb R^n \to \overline{\mathbb R}$ is lower semicontinuous.  This type of problems appear frequently in several practical areas like, e.g., compressed sensing \cite{yang2011alternating}, matrix completion \cite{6262492},  signal processing \cite{Chartrand2007,BrucksteinDonohoElad2009}, and many more.

A  typical model function is
\begin{equation}\label{Eq:ModelforAC}
f_{\bar x}(x):=h(\bar x)+ \left<\nabla h(\bar x), x-\bar x \right>+q(x),
\end{equation}
 for which the local error model reduces to 
\begin{equation}\label{Eq:MAP}
\left|f(x)-f_{\bar x}(x)\right| =\left|h(x)-h(\bar x)-\left<\nabla h(\bar x), x-\bar x \right>\right|,
\end{equation}
i.e., it depends on the degree of smoothness of $h$ only.  Let us consider for all $x \in B_{a}(\bar x)$ with some constant $a>0$,  in particular,  when $\nabla h$ is $\psi$-uniformly continuous on $B_{a}(\bar x)$, the local error can be bounded by $\int_0^1 \psi(s|x-\bar x|)|x-\bar x| \text{ds}$ by \cite[Lemma~3.1]{ochs2019non}, which degenerates into $\frac{L_a}{2}\|x-\bar x\|^2$ if $\nabla h$ is $L_a$-Lipschitz continuous on $B_{a}(\bar x)$.  Note that the above error bound is a relaxation of the global version mentioned in \cite[Example~5.1]{ochs2019non}.  Additionally,  if $\nabla h$ is Lipschitz continuous, then both \cref{Prop:App1} and \cref{Prof:App2} imply that \Cref{Ass:H} is valid. 

If we assume that $\nabla h$ is locally Lipschitz continuous on its domain and the Hessian approximation $H_k:=\identity$, then \Cref{Alg:ModelQuasiNewton} becomes forward-backward splitting or proximal gradient methods,  and our previously obtained results align with those presented in \cite[Section~3]{KanzowMehlitz2022} about the monotone proximal gradient method and \cite{jia2023convergence}.  


\subsection{Composite problems}\label{Subs:CP}
More generally, we consider the following nonconvex nonsmooth composite problems
\begin{equation}\label{Eq:Comp}
	\min_x \ f(x) :=q(x)+h\big(A(x)\big),
\end{equation}
where $q: \mathbb R^n \to \overline{\mathbb R}$ is proper, lower semicontinuous and $h: \mathbb R^m \to \mathbb R$ is continuously differentiable, and $A:\mathbb R^n \to \mathbb R^m$ is a possibly nonlinear $\mathcal C^1$ mapping over $\mathbb R^n$.  The notable examples include low-rank matrix recovery problems \cite{gillis2014and,charisopoulos2021low,kaplan2021low, garber2023linear}, quadratic inverse problems \cite{duchi2018stochastic,BolteSabachTeboulleVaisbourd2018,flemming2018variational},  image processing \cite{BianChen2015}, and so on. 
\subsubsection{Model function w.r.t.  the linearization of $A$}\label{Ssubs:A}
 Let us consider the (linear) Taylor approximation of $A$, then the problem \eqref{Eq:Comp} can be modeled as
\begin{equation}\label{Eq:ModelEX2.1}
f_{\bar x}(x):=h\big(A(\bar x)+\mathcal D A(\bar x)(x-\bar x)\big)+q(x),
\end{equation}
where the local model error reduces to 
\begin{equation}\label{Eq:MAP2}
|f(x)-f_{\bar x}(x)| =\big|h\big(A(x)\big)-h\big(A(\bar x)+\mathcal D A(\bar x)(x-\bar x)\big)\big|,
\end{equation}
where $x \in B_{a}(\bar x)$ with some constant $a>0$.  When $ h$ is $L$-Lipschitz continuous, the error can be bounded by $L|A(x)-A(\bar x)-\mathcal D A(\bar x)(x-\bar x)|$,  also if $\mathcal D{A}$ is $\beta_a \text{-Lipschitz}$ on $B_{a}(\bar x)$,  the error can be bounded by $\frac{L\beta_a}{2}\|x-\bar x\|^2$.  Since $A$ is continuously differentiable, for $x$ sufficiently close to $\bar x$, both $A(x)$ and $A(\bar x)+\mathcal D A(\bar x)(x-\bar x)$ lie in a neighborhood of $A(\bar x)$, in this case,  the local Lipschitz continuity of $h$ is also valid \cite[Example~5.2]{ochs2019non}. 
By \cref{Prof:App2},  \Cref{Ass:H} holds provided that $\nabla h$ and $\nabla (h\circ A)$ are Lipschitzly continuous.

Our method is related to the proximal decent methods \cite{drusvyatskiy2018error,lewis2016proximal}.  Note that if $q$ and $h$ are convex,  where $h$ is L-Lipschitz continuous and the Jacobian $\mathcal D A$ is $\beta_{a}$-Lipschitz continuous,  and $H_k:=\identity$, then \Cref{Alg:ModelQuasiNewton} becomes the proximal descent methods as in \cite[Algorithm~1]{drusvyatskiy2018error} and our convergence results in this situation can cover the most results in \cite[Section~5]{drusvyatskiy2018error}.  If $h$ is just continuous, not convex anymore, and $A$ is not necessarily assumed to have a local Lipschitz gradient, then the convergence and rate-of-convergence can be obtained by the aforementioned sections, which, however, were not discussed in \cite{drusvyatskiy2018error}. 

More generally, if we assume that $h$ is not differentiable and $q:=0$, then \eqref{Eq:Comp} becomes the problem discussed in \cite{lewis2016proximal} and \cite[Section~8]{drusvyatskiy2018error}.  However, in \cite[Section~8]{drusvyatskiy2018error}, just the rationality of the corresponding linear convergence of the algorithms was predicted, where no convincing proof was provided.  In \cite{lewis2016proximal}, just the global convergence of the algorithm was obtained, not the rate-of-convergence included. 
\subsubsection{Model function w.r.t. the linearization of $h$}\label{Subsec:Ex2.2}
In this case, we can define the model function of \eqref{Eq:Comp} at the model center $\bar x$ as
\begin{equation}\label{Eq:modelEx2}
f_{\bar x}(x):=h\big(A(\bar x)\big)+\big(\nabla h \big(A(\bar x)\big)\big)^T\big(A(x)-A(\bar x)\big)+q(x),
\end{equation}
for which the local error model can be formulated as 
\begin{equation}\label{Eq:MAP3}
|f(x)-f_{\bar x}(x)| =\big|h\big(A(x)\big)-h\big(A( \bar x)\big)-\big(\nabla h \big(A(\bar x)\big)\big)^T\big(A(x)-A(\bar x)\big)\big|\\
\end{equation}
where $x \in B_{a}(\bar x)$. When $ \nabla h$ is {L}-Lipschitz continuous, the error can be bounded by $\frac{L}{2}\|A(x)-A(\bar x)\|^2$,  also if $A$ is $\beta_a$-Lipschitz continuous on $B_a(\bar x)$, the error can be further bounded by $\frac{L\beta_a^2}{2}\|x-\bar x\|^2$. This is motivated by \cite{ochs2015iteratively}, however, where requires that the part $f_{\bar x}(x)-q(x)$ is convex and nondecreasing to guarantee the existence of the solution(s) of the corresponding subproblems, which is not necessary in our manuscript since the subproblems \eqref{Eq:NonSubki} of \Cref{Alg:ModelQuasiNewton} include the proximal part as well as uniform positive definiteness of $H_k$ in view of \Cref{Ass:ModelNewton}\ \ref{item: Unipd}, both together for the existence of the solutions of the subproblems.  
By \cref{Prop:App1},  \Cref{Ass:H} holds provided that $\nabla h$ and $A$ are Lipschitz continuous.

Our method can be degenerated to the IRLS algorithm \cite[Algorithm~6]{ochs2015iteratively} (when $H_k:=\identity$), where the convergence was illustrated in \cite[Theorem~2]{ochs2015iteratively} by assuming that $h$ has locally Lipschitz gradients, $g$ is convex,  and $f$ has the \rm{KL} property.  So we definitely say that our desired convergence results cover the convergence analysis in \cite[Section~5]{ochs2015iteratively}. In addition, if $A(x):=Mx-b$ with a matrix $M \in \mathbb R^{m \times n}$ and a vector $b \in \mathbb R^m$,  \eqref{Eq:Comp} covers the optimization problem discussed in \cite{liu2024inexact}, it requires $g$ is convex, $h$ is twice continuously differentiable on an open set containing $M(\mathcal O)-b$ where $\mathcal O$ is an open set covering $\dom g$,  and $f$ is coercive.  Personally, those assumptions are too restricted, we do not require any in this manuscript. In \cite{liu2024inexact}, the Hessian approximation has been required to be uniformly bounded \cite[Lemma~4]{liu2024inexact}, which is necessary for using the KL property \cite[Theorem~4]{liu2024inexact}.  However, Our work do not require any boundedness of the Hessian approximation.

In some practical areas, like signal processing, 
 machine learning,
\ compressed sensing, 
and image processing, 
typically,  $g$ is the regularization function used to promote the sparser structure of the solution(s),  and $h$ is always non-negative, which motivates us to introduce the following two model functions,
\begin{equation*}
f_{\bar x}(x):=\max\big\{0,h\big(A(\bar x)\big)+\big(\nabla h \big(A(\bar x)\big)\big)^T\big(A(x)-A(\bar x)\big)\big\}+q(x),
\end{equation*}
and 
\begin{equation}\label{Eq:ModelExabs}
f_{\bar x}(x):=\big|h\big(A(\bar x)\big)+\big(\nabla h \big(A(\bar x)\big)\big)^T\big(A(x)-A(\bar x)\big)\big|+q(x).
\end{equation}
Both reserve the nonnegative property of $h$. 
\section{Numerical Experiments}
We implemented \Cref{Alg:ModelQuasiNewton},  based on the underlying choice of the model function proposed in \Cref{Sec:Ex}, in MATLAB (R2023b) and tested it on some practical problems.    All test runs use the following parameters
\begin{equation*}
\tau:=2, \ \delta:=0.25, \ \mu:=0.5, \ \gamma_{\min}:=1, \ \gamma_{\max}:=10^{10}.
\end{equation*}
\subsection{Polytope feasibility} 
Polytope feasibility problems aim to find a feasible point $x^* \in \mathcal F$, where $\mathcal F$ is defined as $$\mathcal F:=\{x\in \mathbb R^n \,|\, \left<a_i,x\right> \leq b_i, 1\leq i \leq m \},$$ 
i.e., a polytope in $\mathbb R^n$, by minimizing the following optimization problem
\begin{equation}\label{Eq:PolytopeFeas}
\min_{x\in \mathbb R^n } f(x):=\sum_{i=1}^m \left(\left<a_i,x\right>-b_i \right)_+^p,
\end{equation}
where $(\cdot)_+:=\max\{0,\cdot\}$ is a positive slicing and $p \geq 2$ is a given parameter.  The function $f$ is known to satisfy the KL property. 

To approximate (even smooth) the nonsmooth term $\left(\left<a_i,x\right>-b_i \right)_+$, we may consider the so-called normalized softplus function $\phi_c(t):=\frac{1}{c}\log\left(1+e^{ct}\right)$, which satisfies $\phi_c(t)\geq (t)_+$ for all $c>0$ and $t\in \mathbb R$. In particular, $$\frac{1}{c}\log\left(1+e^{c\left(\left<a_i,x\right>-b_i\right)}\right) \geq \left(\left<a_i,x\right>-b_i \right)_+ \quad \forall 1 \leq i \leq m.$$Therefore it is natural  to consider the first-order Taylor expansion of the function $\sum_{i=1}^m \phi_c^p\left(\left<a_i,x\right>-b_i\right)$ as a model function of $f$.  More precisely,  for each point $\bar x\in \mathbb R^n$, the model error can be explicitly calculated as
\begin{equation*}
\begin{aligned}
&~~~~\left | \sum_{i=1}^m \left(\left<a_i,x\right>-b_i \right)_+^p-c^{-p}\sum_{i=1}^m\log^p\left(1+e^{c\left(\left<a_i,\bar x\right>-b_i\right)}\right)-\left<\text{grad},x-\bar x\right>\right|\\
&\leq \left|c^{-p}\sum_{i=1}^m\log^p\left(1+e^{c\left(\left<a_i, x\right>-b_i\right)}\right)-c^{-p}\sum_{i=1}^m\log^p\left(1+e^{c\left(\left<a_i,\bar x\right>-b_i\right)}\right)-\left<\text{grad},x-\bar x\right>\right|\\
& \leq \frac{L_{\bar x}}{2}\|x-\bar x\|^2 \quad \forall x \ \text{around}\ \bar x,
\end{aligned}
\end{equation*}
where $\text{grad}:=c^{1-p} p\sum_{i=1}^m \frac{e^{c\left(\left<a_i,\bar x\right>-b_i\right)}}{1+e^{c\left(\left<a_i,\bar x\right>-b_i\right)}}\log^{p-1}\left(1+e^{c\left(\left<a_i,\bar x\right>-b_i\right)}\right)a_i$ and $L_{\bar x}$ is the Lipschitz constant of $\text{grad}$ at $\bar x$.

We compare the performance of \Cref{Alg:ModelQuasiNewton} with that of the corresponding gradient method.  In \Cref{Alg:ModelQuasiNewton}, the Hessian approximation in \Cref{Alg:ModelQuasiNewton} is updated via
\begin{equation*}
H(x^k):=P_{\mathcal S_{+}(\mathbb R^2)}(\hat H_k)+\mu \identity,
\end{equation*}
where $\hat H_k$ denotes the Hessian matrix of the function $\sum_{i=1}^m \phi_c^p\left(\left<a_i,x\right>-b_i\right)$ at the iteration $k$.  For the gradient method, we simply choose $h(x^k):=\text{Id}$ for every iteration.

For implementation, we set $n=100$, $m=200$. The data $a_i$ and $b_i$ is sampled independently from random uniform distribution on $[-1,1]$ for all $1\leq i\leq m$.  The starting point is fixed as $x_0=(1,\ldots, 1)^T\in \mathbb R^n$.  All settings are the same as those in \cite[Section 7]{doikov2024super}.  The parameter $c$ for the normalized softplus function is chosen to be $c=2$.  We tested the two types of algorithms on polytope problems with various values of $p\in\{2,3,3.5,4\}$.  The termination criterion is set as $\sum_{i=1}^m \left(\left<a_i, x\right>-b_i\right) \leq 1e-4$.  Notably, the gradient method encountered numerical issues and returned ``NaN'' errors across all tested $p$. As a result, only the results of \Cref{Alg:ModelQuasiNewton} are reported (zoomed in), which is shown below.
\begin{figure}[H]
  \centering
    \includegraphics[width=0.65\textwidth]{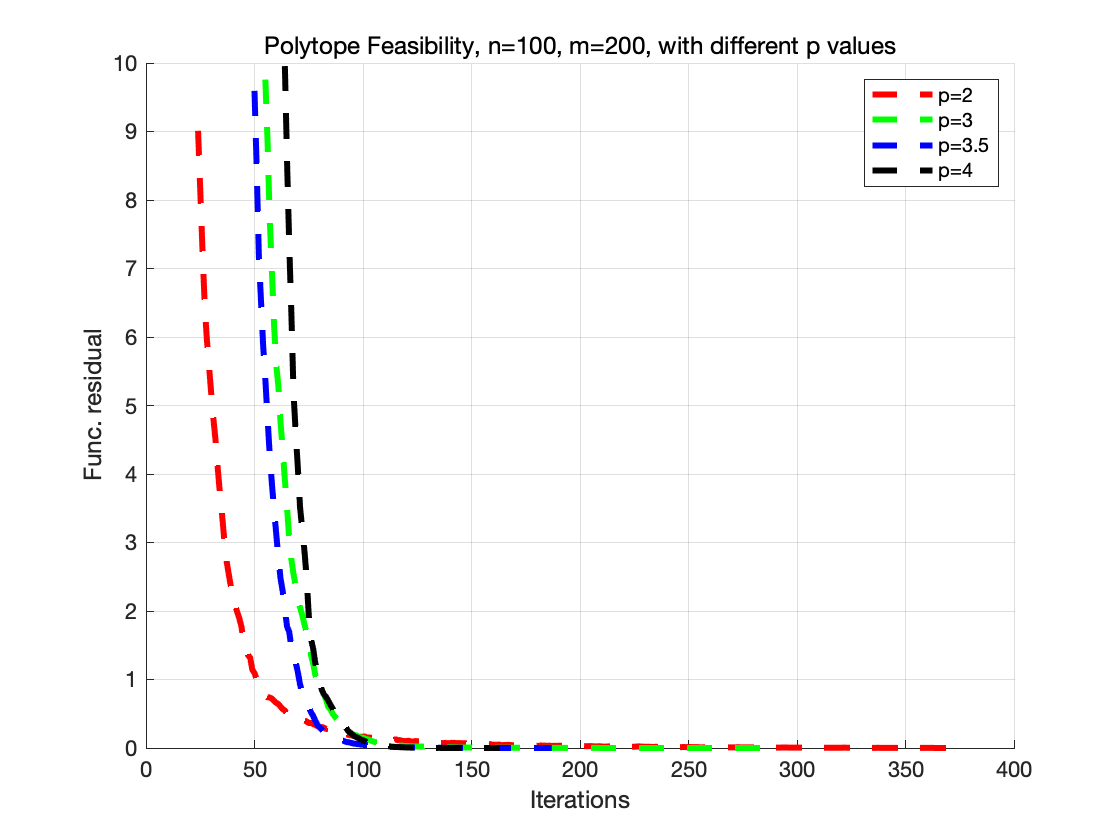}
\end{figure}
The numerical experiments demonstrate that \Cref{Alg:ModelQuasiNewton} consistently drives the objective function value close to zero within 100 iterations for each tested instance. Furthermore, the results indicate that as $p$ increases, the objective function becomes steeper near the solution, leading to faster convergence of the algorithm.
\subsection{Quadratic inverse problems}
Quadratic inverse problems aim essentially to solve approximately a system of quadratic equations \cite{BolteSabachTeboulleVaisbourd2018, mukkamala2022global,censor1981iterative,luke2017phase,wang2017solving}.  Let the so-called sampling matrix $A_i \in \mathbb R ^{n \times n}, i=1,\ldots,m$ be symmetric positive semi-definite and possibly noisy measurements $b_i\in \mathbb R^n$,  the goal of quadratic inverse problem is to find $x\in \mathbb R^n$ satisfying $x^TA_ix \approx b_i$ for all $i=1,\ldots, m$.  Adopting a quadratic function to measure the error, the problem can then be reformulated as the following nonconvex and nonsmooth optimization problem
\begin{equation}\label{Eq:PR}
\min_{x\in \mathbb R^n} h(A(x)):=\frac{1}{m} \sum_{i=1}^m \frac{1}{2}\left(x^TA_ix-b_i\right)^2.
\end{equation}
Meanwhile, we also consider the corresponding sparse quadratic inverse problem which aims to obtain the sparser solution
\begin{equation}\label{Eq:SPR}
\min_{x\in \mathbb R^n} \frac{1}{m} \sum_{i=1}^m \frac{1}{2}\big(x^TA_ix-b_i\big)^2+\lambda \|x\|_1,
\end{equation}
where $\lambda>0$ plays the role of a penalty parameter controlling the trade-off between the system fidelity versus its regularizer $\|\cdot\|_1$, and $\|\cdot\|_1$ is the number of nonzero element. Note that objective functions of both problems \eqref{Eq:PR} and \eqref{Eq:SPR} can be bounded below by $0$.

Clearly, the analysis falls under the category of composite problems \eqref{Eq:Comp}, where $h(A(x)):=1/{(2m)} \sum_{i=1}^m (x^TA_ix-b_i)^{2}$ and $g(x):=\lambda \|x\|_1$. 
We consider the following three model functions to solve the problems in \eqref{Eq:SPR}.

\textbf{Model 1 (M1).} As mentioned above, the problem falls under the structure of additive composite problems in \Cref{Subs:CP}. 
 Consider the standard model function \eqref{Eq:modelEx2} for the problem \eqref{Eq:SPR}, where around $x^k \in \mathbb R^n$,  the model function $f^1: \mathbb R^n \to \mathbb R$ is given by
\begin{equation*}
f_{x^k}^1(x):=\frac{1}{2m} \sum_{i=1}^{m} \left({x^k}^TA_i{x^k}-b_i\right)^{2 }+2\left({x^k}^TA_ix^k-b_i\right) \left<A_ix^k, x-x^k \right>+\lambda \|x\|_1,
\end{equation*}
which is the approximate first-order Taylor expansion with respect to the first term in \eqref{Eq:SPR}.  Note that this model function is always lower semicontinuous.  

\textbf{Model 2 (M2)}. The importance of finding better model functions is to make the model function closed to the actual objective function; the closer is, the better is.  Naturally,  we consider the associated second-order approximation expansion (without second-order Hessian information), and the model function $f^2: \mathbb R^n \to \mathbb R$, centered at $x^k$,  is given by
\begin{equation*}
\begin{aligned}
 f_{x^k}^2(x)&:=\frac{1}{m} \sum_{i=1}^{m} \frac{1}{2}\left({x^k}^TA_i{x^k}-b_i\right)^{2 }+2\left({x^k}^TA_ix^k-b_i\right) \left<A_ix^k, x-x^k \right>\\
&~~~~+\frac{1}{2}(x-x^k)^T(x-x^k)+\lambda \|x\|_1.
\end{aligned}
\end{equation*}
The model function is still lower semicontinuous. 

{\textbf{Model 3 (M3).} The another principle for better choice of the model function is to reserve the property of the original function as much as one can.  For the problem \eqref{Eq:SPR}, $h(A(x)):=1/{(2m)} \sum_{i=1}^m (x^TA_ix-b_i)^{2}$ is always nonnegative, so we consider the following model function $f^3:\mathbb R^n \to \mathbb R$ from \eqref{Eq:ModelExabs}, centered at $x^k$,
\begin{equation*}
f_{x^k}^3(x):=\frac{1}{2m} \sum_{i=1}^{m} \left|\left({x^k}^TA_i{x^k}-b_i\right)^{2 }+2\left({x^k}^TA_ix^k-b_i\right) \left<A_ix^k, x-x^k \right>\right|+\lambda \|x\|_1.
\end{equation*}
This function is also lower semicontinuous,  we solve the subproblems by the primal-dual hybrid gradient method (PDHG) \cite{6126441}, where parameters are chosen the same as those in \cite[Section~5.2]{mukkamala2022global}.

Note that all the above model functions (M1, M2, M3) are obviously convex, hence \cref{Ass:ModelNewton} \ref{item: mf} is valid.  Through easy calculations,  we can verify that for M1, M2,  M3 (an additional termination criterion is employed),  \Cref{Ass:H} holds.
Let us consider the principle operator of $H$ in \Cref{Alg:ModelQuasiNewton}, for comparison, we here choose the following two update principles.  The first one is the Hessian approximation, i.e.,
\begin{equation}\label{Eq:HA}
H(x^k):=P_{\mathcal S_{+}(\mathbb R^n)}(\nabla^2 h(A(x^k)))+\mu \identity,
\end{equation}
for all $k=1,\ldots$, where $h(A(x^k)):=1/{(2m)} \sum_{i=1}^m ((x^k)^TA_ix^k-b_i)^{2}$.  In this case, we call \Cref{Alg:ModelQuasiNewton} as model proximal quasi-Newton methods (MQN for short).  In this case, we employ an alternating direction method of multipliers (ADMM) method to solve the subproblem for M1 and M2.

The next one is
\begin{equation*}
H_{k}=L_{k}\identity,
\end{equation*}
for all $k=1, \ldots$,  where $L_{k}$ is the Barzilai and Borwein stepsize \cite{barzilai1988two}, then \Cref{Alg:ModelQuasiNewton} degenerates into the so-called model proximal gradient methods (MG for short).  In this case, when we use M1 and M2 as the model function, solutions of the subproblems can be calculated by the so-called soft-thresholding operator. 

We used 100 random synthetic datasets where $n:=50$ and $m:=1000$ to test \Cref{Alg:ModelQuasiNewton} and compare the empirical results generated by the model quasi-Newton method ($H_{k+1}$ is updated by \eqref{Eq:HA}) and the model gradient method ($H_{k+1}=L_{k+1}\identity$) by employing the different model functions (Model 1, Model 2,  Model 3), respectively.  The initial stepsize for each iteration is selected as $\gamma_k^0:=2$ for all $k \geq 1$ and $\gamma_0^0:=\|\nabla h(A(x_0))\|_{\infty}, \|\nabla h(A(x_0))\|_{2}, \|\nabla h(A(x_0))\|_{2}$ for Model 1, Model 2,  Model 3, respectively.  We terminated the algorithms where M1 or M2 is employed as the model function, if
\begin{equation*}
\frac{f(x^k)-f(x^{k+1})}{\max\{1, f(x^{k+1})\}} \leq 10^{-4}.
\end{equation*}
It is together with
$$\sum_{i=1}^{m}\left|\left({x^k}^TA_i{x^k}-b_i\right)^{2 }+2\left({x^k}^TA_ix^k-b_i\right) \left<A_ix^k, x^{k+1}-x^k \right>\right|\leq 10^{-4}$$
as the termination criterion when M3 is employed.

Using the vector $0.1*\text{ones} (n,1)$ as the starting point and testing three values $\{e^{-2},e^{-3},e^{-4}\}$ for parameter $\lambda$ for all 100 testproblems,  we reported the average results shown in \Cref{Table:PR}.  The average number of (outer) iterations is denoted by $k$, $j$ is the average accumulated number of the backtracking line search,  CPU means the average total cost time in seconds, $f_v$ denotes the optimal functional value on average,   $d_f$ denotes the average model error at the last iteration (the distance between objective function value and its model function value at the last iteration), and $r$ denotes of the average numbers of nonzero components of the obtained solutions.  $``-"$ means that the corresponding algorithm can not converge in 2000 iterations for at least one testproblem.\\
\begin{scriptsize}
	\def\tablename{\normalsize Table}
\begin{table}[htbp]\caption{Averaged numerical results for 100 random quadratic inverse problems.}\label{Table:PR}
\begin{center}
\begin{tabular}[H]{r|c|rrrrrr}
\hline
      $\lambda$ & Algorithm& $ k $ & $ j $ & CPU(s) & $f_v\hspace{0.6mm}$&  $d_f\hspace{0.6mm} $& r \\ \hline
\multirow{6}{*}{$1e^{-2}$} & MQN-M1 & 2 &    5.08 &    0.2745 &      0.0049 &     0 &   0\\
&  MQN-M2 & 3 &    6 &   0.3308&    0.0050 &     0 &    0 \\
&  MQN-M3& 6.05&    7.23&     1.1755&     0.0048&     4.85e-8&     2\\
\cline{2-8}
&   MG-M1  & 36.31&  10.18&   3.8759&   0.0051&   0 &   0\\
& MG-M2  & 25.61&   45.06&  3.1124&  0.0049&  1.31e-7& 0\\ 
& MG-M3  &  7.74&   6.49&   0.9266& 0.0062&  1.11e-7& 6.49\\ 
\hline
\multirow{6}{*}{$1e^{-3}$} & MQN-M1 & 2.9&  7.83& 0.3522&  0.0050 &    2.44e-8 &   0 \\
& MQN-M2 & 3&  6&   0.3345&  0.0049&    0 &    0\\
& MQN-M3 & 7.51& 11.01&    1.6256&    0.0055&   8.13e-6 &    16.03\\
\cline{2-8}
& MG-M1 & $-$&  $-$&  $ -$& $  -$& $-$& $ -$\\
& MG-M2 & 48.54 & 91.25&  6.9404&   0.0052&  3.25e-6 &  47.58\\
& MG-M3  &10.09& 14.74& 1.4388& 0.0058& 5.70e-5& 22\\
\hline
\multirow{6}{*}{$1e^{-4}$}& MQN-M1 & 2&  5.11&  0.2618&  0.0051&  0 &    0 \\
&MQN-M2 & 3&  6&   0.3407&  0.0049&  0&    0 \\
& MQN-M3 & 16.50&  31.97&  4.8319&  0.0055&  2.07e-5 &    46.44 \\
\cline{2-8}
& MG-M1 & $-$&  $-$&  $ -$& $  -$& $-$& $ -$\\
& MG-M2& 53.98&   103.37&  8.0046& 0.0053&   7.21e-6& 49.88 \\
& MG-M3& 30.34& 60.84&  6.2678&  0.0056&  9.46e-5&  46.71 \\
\hline
\end{tabular}
\end{center}
\end{table}
\end{scriptsize}
\indent
\Cref{Table:PR} illustrates that when M1 is selected, MG is easier to fail to converge for the smaller $\lambda$, M2 usually needs more inner iterations for both MQN and MG.  Generally,  the corresponding model error were sufficiently small when the algorithms terminate with convergence.  MQN requires fewer outer iterations and fewer inner iterations on average for each outer iteration than MG,  and it generates sparser solutions than MG.

\section{Conclusion}
We presented the so-called model quasi-Newton method for addressing the nonconvex and nonsmooth optimization problems. This algorithm combines the proximal minimization of the (local) model function and the backtracking line search to ensure that the sequential decrease of the objective function.  In this manuscript, we did not impose any assumption of boundedness on the sequence of variable metrics since it is too restrictive even illogical,  particularly for the objective function with sharp curves, instead, we required variable metrics are generated by a continuous matrix generator.  By assuming the first-order information of the model function,  we obtained the subsequential convergence, where the sequence of variable metrics is not uniformly bounded.  Furthermore, by employing the KL property at the accumulation point of the generated iterative sequence,  the convergence of the entire sequence to a stationary point and the corresponding rate-of-convergence were established.  Those illustrated that the boundedness Hessian approximation should be a problem-tailored consequence of convergence results,  which is not logical to be assumed as a prerequisite for the convergence analysis. We also provided examples of the local model functions for different types of (addictive) composite problems to empower the generality of our optimization problem and algorithm. Numerically, we compared our algorithm with its associate gradient method to tackle the polytope feasibility problems and (sparse) quadratic inverse problems (employing the different classes of model functions),  both problems illustrated the effectiveness and robustness of our algorithm. 

In the future, on the one hand, we will focus on the theoretical analysis that the boundedness of the sequence of variable metrics (or Hessian approximation) is the byproduct of problem-tailored convergence results, particularly, dependent on the regularity of the objective function and the specific (proximal) quasi-Newton methods.  On the other hand, in the view of learning to optimize (L2O), we focus on learning the Hessian approximation based on prior information for the proximal quasi-Newton methods in order to improve their overall effectiveness.

\smallsection{Acknowledgments} {\small
Xiaoxi Jia would like to thank Shida Wang for his valuable discussions on the (quasi-)Newton-type methods.
}

\section{Appendix}\label{Appendix}
Let us consider 
\begin{equation}\label{Eq:cf}
f(x):=q(x)+h\big(A(x)\big),
\end{equation}
where $q: \mathbb R^n \to \overline{\mathbb R}$ is proper, lower semicontinuous and $h: \mathbb R^m \to \mathbb R$ is continuously differentiable, and $A:\mathbb R^n \to \mathbb R^m$ is a possibly nonlinear $\mathcal C^1$ mapping over $\mathbb R^n$.  We next give some special examples of model functions under the mild requirements to make sure \eqref{H2} holds for all $\bar x \in \dom f$.
\begin{proposition}\label{Prop:App1}
Let $f$ be defined as \eqref{Eq:cf}. For any fixed $\bar x\in \dom f$, suppose the corresponding model function is given by
\begin{equation*}
f_{\bar x}(x):=h\big(A(\bar x)\big)+\left\langle \nabla h\big(A(\bar x)\big), A(x)- A(\bar x) \right\rangle+q(x) \quad \forall x\in \dom f.
\end{equation*}
Then \eqref{H2} holds if $\nabla h$ and $A$ are locally Lipschitz continuous.
\end{proposition}
\begin{proof}
We now have
\begin{equation*}
g_{\bar x}(x)=h\big(A(\bar x)\big)+\left\langle \nabla h\big(A(\bar x)\big), A(x)-A(\bar x) \right\rangle-h\big(A(x)\big),
\end{equation*}
which is smooth in terms of $x$.  Then,
\begin{equation*}
\begin{aligned}
\nabla g_{\bar x}(x)&=\left(\mathcal D A(x)\right)' \nabla h\big(A(\bar x)\big)-\left(\mathcal D A(x)\right)' \nabla h\big(A(x)\big)\\
&=\left(\mathcal D A(x)\right)'\left(\nabla h\big(A(\bar x)\big)-\nabla h\big(A(x)\big)\right).
\end{aligned}
\end{equation*}
Assuming $L_1>0$ is the Lipschitz constant of $\nabla h$ and $L_2>0$ is the Lipschitz constant of $A$, then $\|\mathcal D A(x)\| \leq L_2$, so
\begin{equation*}
\|\nabla g_{\bar x}(x)\| \leq L_1(L_2)^2\|x-\bar x\| \quad \forall x\ \text{closed to}\ \bar x,
\end{equation*}
this shows that condition \eqref{H2} holds with $L:=L_1(L_2)^2$.
\end{proof}
\begin{proposition}\label{Prof:App2}
Let $f$ be defined as \eqref{Eq:cf}. For any fixed $\bar x\in \dom f$, suppose the corresponding model function is given by
\begin{equation*}
f_{\bar x}(x):=h\big(A(\bar x)+ \mathcal DA(\bar x)(x-\bar x) \big)+q(x) \quad \forall x\in \dom f.
\end{equation*}
Then \eqref{H2} holds if $\nabla h$ and $\nabla (h\circ A)$ are locally Lipschitz continuous.
\end{proposition}
\begin{proof}
For any fixed $\bar x\in \dom f$, we have
\begin{equation*}
g_{\bar x}(x)=h\big(A(\bar x)+\mathcal D A(\bar x)(x-\bar x) \big)-h\big(A(x)\big) \quad \forall x\in \dom f,
\end{equation*}
which is differentiable in terms of $x$.  Then
\begin{equation*}
\begin{aligned}
\nabla g_{\bar x}(x)&=\left(\mathcal D A(\bar x)\right)' \nabla h\big(A(\bar x)+ \mathcal D A(\bar x)(x-\bar x) \big)-\left(\mathcal D A( x)\right)'\nabla h\big(A(x)\big),\\
&=\left(\mathcal D A(\bar x)\right)'\left( \nabla h\big(A(\bar x)+ \mathcal D A(\bar x)(x-\bar x)\big) -\nabla h \big(A(\bar x)\big)\right)\\
&~~~~+\left(\mathcal D A(\bar x)\right)' \nabla h\big(A(\bar x)\big)-\left(\mathcal D A( x)\right)'\nabla h\big(A(x)\big),\\
\end{aligned}
\end{equation*}
By assuming $L_1>0$ is the Lipschitz constant of $\nabla h$,  $L_2>0$ is the Lipschitz constant of $\nabla (h\circ A)$, we have
\begin{equation*}
\begin{aligned}
\|\nabla g_{\bar x}(x)\| &\leq \left\|\left(\mathcal D A(\bar x)\right)'\left( \nabla h\big(A(\bar x)+ \mathcal D A(\bar x)(x-\bar x)\big) -\nabla h \big(A(\bar x)\big)\right)\right\|\\
&~~~~+\left\|\left(\mathcal D A(\bar x)\right)' \nabla h\big(A(\bar x)\big)-\left(\mathcal D A( x)\right)'\nabla h\big(A(x)\big)\right\|\\
&\leq L_1\|\mathcal D A(\bar x)\|^2\|\|x-\bar x\|+L_2\|x-\bar x\|\\
&=\left(L_1\|\mathcal D A(\bar x)\|^2+L_2\right)\|x-\bar x\| \quad \forall x \ \text{closed to}\ \bar x.
\end{aligned}
\end{equation*}
Therefore, the condition \eqref{H2} holds with $L:=L_1\|\mathcal D A(\bar x)\|^2+L_2$.
\end{proof}

Actually, if $A$ is the identity operator, then \eqref{H2} holds only under $\nabla h$ is Lipschitz continuous. 
\bibliographystyle{habbrv}
\bibliography{references}

\end{document}